\documentclass[1 [leqno,11pt]{amsart}
\usepackage{amssymb, amsmath,latexsym,amsfonts,amsbsy, amsthm,mathtools,graphicx,CJKutf8,CJKnumb,CJKulem,color,xcolor,mathrsfs,dsfont}
\usepackage{float}
\usepackage{hyperref}
\usepackage{braket}
\usepackage{graphicx}
\usepackage{bbm}
\renewcommand{\theequation}{\thesection.\arabic{equation}}
\numberwithin{equation}{section}

\allowdisplaybreaks

\let\al=\alpha

\let\d=\delta

\let\la=\lambda

\let\f=\frac

\let\om=\omega

\let\La=\Lambda

\let\pa=\partial

\newcommand{\cal}{\mathcal}
\def\bH{\boldsymbol{H}}
\def\bV{\boldsymbol{V}}

\def\cL{{\cal L}}
\def\cM{{\mathcal M}}
\def\cN{{\mathcal N}}
\def\cO{{\cal O}}
\def\cP{{\mathcal P}}

\def\cR{{\mathcal R}}
\def\cS{{\mathcal S}}

\def\R{\mathbb R}
\def\Z{\mathbb Z}
\def\T{\mathbb T}
\def\N{\mathbb N}
\def\G{{\scriptscriptstyle G}}

\def\dive{\mathop{\rm div}\nolimits}

\def\grad{\nabla}
\def\lap{\Delta}
\def\lapxi{\lap_{ \scalebox{0.8}{\xxi}}}
\def\gradxi{\nabla_{\scalebox{0.8}{\xxi}}}
\def\sxxi{{\scalebox{0.6}{\xxi}}}

\def\oo{\infty}
\def\ooplus{\mathop{\scalebox{1.5}{$\oplus$}}}
\def\ccap{\mathop{\scalebox{1.5}{$\cap$}}}

\renewcommand{\d}{\,{\rm d}}            

\newcommand{\uu}{{\boldsymbol{u}}}
\newcommand{\vv}{{\boldsymbol{v}}}
\newcommand{\xx}{{\boldsymbol{x}}}

\newcommand{\eeta}{{\mbox{\boldmath$\eta$}}}
\newcommand{\xxi}{{\mbox{\boldmath$\boldsymbol{\xi}$}}}
\newcommand{\hh}{{\boldsymbol{H}}}

\newcommand{\oome}{{\boldsymbol{\Omega}}}
\newcommand{\bdive}{\overline{\dive}}
\newcommand{\blap}{\overline{\lap}}
\newcommand{\bgrad}{\overline{\grad}}
\newcommand{\UU}{\boldsymbol{U}}
\newcommand{\sL}{\mathscr{L}}
\newcommand{\ltr}{L^2_r}
\def\eqdefa{\buildrel\hbox{\footnotesize def}\over =}
\newcommand{\vg}{\vv^\G}

\let\f=\frac

\def\Re{{\rm Re\,}}
\def\Im{{\rm Im\,}}

\def\sign{\mathop{\rm sign}}

\def\epsilon{\varepsilon}
\def\phi{\varphi}

\newcommand{\tr}{{\rm T}}

\newcommand{\andf}{\quad\hbox{and}\quad}
\newcommand{\with}{\quad\hbox{with}\quad}

\newcommand{\beq}{\begin{equation}}
\newcommand{\eeq}{\end{equation}}
\newcommand{\ben}{\begin{eqnarray}}
\newcommand{\een}{\end{eqnarray}}
\newcommand{\beno}{\begin{eqnarray*}}
\newcommand{\eeno}{\end{eqnarray*}}

\usepackage{float}

\newtheorem{Theorem}{Theorem}[section]

\newtheorem{Proposition}[Theorem]{Proposition}
\newtheorem{Lemma}[Theorem]{Lemma}

\begin{document}
\title[Nonlinear optimal asymptotic decay around the 3D Oseen vortex filament]
{Nonlinear asymptotic stability and optimal decay rate around the three-dimensional Oseen vortex filament}

\author{Te Li, Ping Zhang, and Yibin Zhang}

\address[T. Li]
 {State Key Laboratory of Mathematical Sciences, Academy of Mathematics  \&  Systems Science,   Chinese Academy of Sciences,  Beijing 100190,China.}
\email{teli@amss.ac.cn}

\address[P. Zhang]
{State Key Laboratory of Mathematical Sciences, Academy of Mathematics  \&  Systems Science,   Chinese Academy of Sciences,  Beijing 100190,China \newline
and \newline
 School of Mathematical Sciences, University of Chinese Academy of Sciences  Beijing 100049,    China}
\email{zp@amss.ac.cn}

\address[Y. Zhang]
 {Academy of Mathematics $\&$ Systems Science, The Chinese Academy of
Sciences, Beijing 100190, CHINA.  } \email{zhangyibin22@mails.ucas.ac.cn}

\date{\today}

\begin{abstract}
In the high-Reynolds-number regime, this work investigates the long-time dynamics of the three-dimensional incompressible Navier–Stokes equations near the Oseen vortex filament. The flow exhibits a strong interplay between vortex stretching, shearing, and mixing, which generates ever-smaller spatial scales and thereby significantly amplifies viscous effects. By adopting an anisotropic self-similar coordinate system adapted to the filament geometry, we establish the nonlinear asymptotic stability of the Oseen vortex filament. All non-axisymmetric perturbations are shown to decay at the optimal rate $t^{-\kappa |\alpha|^{1/2}}$. At the linear level, this decay mechanism corresponds to a sharp spectral lower bound $\Sigma(\alpha) \sim |\alpha|^{1/2}$ for the nonlocal Oseen operator $L_\perp - \alpha \Lambda_\perp$, and we identify an explicit spectral point attaining this optimal bound. Combined with the spectral estimates obtained in \cite{LWZ}, our analysis fully resolves the conjecture proposed in \cite{GM} concerning the asymptotic scaling laws for the spectral and pseudospectral bounds $\Sigma(\alpha)$ and $\Psi(\alpha)$. These results provide a rigorous mathematical explanation for the shear–mixing mechanism in the vicinity of the 3D Oseen vortex filament.

\end{abstract}

\maketitle


\section{Introduction}\label{intro}
In this paper, we consider the three-dimensional (3D) incompressible Navier–Stokes equations in vorticity form:
\begin{equation}\label{eq 1.1}
\left\{
    \begin{aligned}
        &\partial_t \boldsymbol{\eta} + \boldsymbol{v} \cdot \nabla \boldsymbol{\eta} - \boldsymbol{\eta} \cdot \nabla \boldsymbol{v} = \nu \Delta \boldsymbol{\eta},
        \quad (t, \boldsymbol{x}, z) \in \mathbb{R}^+ \times \mathbb{R}^2 \times \mathbb{T}, \\[4pt]
        &\boldsymbol{v} = \nabla \times (-\Delta)^{-1} \boldsymbol{\eta},
    \end{aligned}
\right.
\end{equation}
where $\boldsymbol{\eta}$ and $\boldsymbol{v}$ denote the vorticity and velocity fields of a viscous fluid flow, respectively, and $\nu > 0$ is the viscosity coefficient. We investigate the nonlinear asymptotic stability and the optimal long-time decay rate of solutions to \eqref{eq 1.1} with a 3D vortex filament as initial data:
\begin{equation}\label{eq 1.1a}
    \boldsymbol{\eta}(t=0) = \alpha\, \delta_{\boldsymbol{x}=\boldsymbol{0}},
\end{equation}
in the high-Reynolds-number regime defined by
\begin{equation}\label{eq: 1.3a}
    \operatorname{Re} \eqdefa |\alpha|/\nu \gg 1.
\end{equation}

Generally, for a smooth oriented curve $\mathcal{C} \subset \mathbb{R}^3$, the measure $\delta_{\mathcal{C}}$ is defined as a divergence-free vector-valued measure satisfying, for any test function $\phi \in \left(C_c^\infty(\mathbb{R}^2\times\mathbb{T})\right)^3$,
\begin{equation*}
    \langle \delta_{\mathcal{C}}, \phi \rangle = \int_{\mathcal{C}} \phi \cdot d\vec{s}.
\end{equation*}

Starting from \eqref{eq: 1.3a}, one can rescale \eqref{eq 1.1}--\eqref{eq 1.1a} so that, without loss of generality, we assume throughout this paper that
\begin{equation}
    \nu = 1 \quad \text{and} \quad \alpha \gg 1.
\end{equation}

It is straightforward to observe that the Navier--Stokes equations \eqref{eq 1.1} admit a family of self-similar solutions with the vortex filament \eqref{eq 1.1a} as initial data, known as the 3D Oseen vortex filament:
\begin{equation}\label{eq 1.2}
   \begin{split}
    \eeta(t,\xx)= \al \eeta^{\G}(t,\xx) \eqdefa \frac{\al}{t} \big(
           0,\,0,\,G\big (\xx/\sqrt{t}\big)
       \big)^\tr, \quad \text{where} \ G(\xxi)\eqdefa (4\pi)^{-1}e^{-|\xxi|^2/4},
   \end{split}     
\end{equation} 
and $\alpha \in \R$ denotes the circulation number. The velocity field associated with this Oseen vortex filament is
\begin{equation}\label{eq 1.3}
\begin{split}
  &\vv(t,\xx) = \al  \boldsymbol{\mathrm{v}}^{\G} (t,\xx)\eqdefa    \frac{\al}{\sqrt{t}} \begin{pmatrix}
      \vg\big(\xx/\sqrt{t}\big)\\0
    \end{pmatrix}, \quad \text{where} \ \vg(\xxi)\eqdefa \frac{\xxi^\perp}{2\pi|\xxi|^2}\big(1-e^{-|\xxi|^2/4}\big).
\end{split}
\end{equation}

Owing to the self-similar structure of the Oseen vortex in the transverse variable $\xx$, it is natural to reformulate \eqref{eq 1.1} in self-similar coordinates
\[
\tau \eqdefa \log t, \qquad \xxi \eqdefa \xx/\sqrt{t},
\]
and to introduce the rescaled vorticity and velocity fields
\[
\bH(\tau,\xxi,z) \eqdefa e^{\tau}\eeta(e^{\tau},e^{\tfrac{\tau}{2}}\xxi,z),
\qquad
\bV(\tau,\xxi,z) \eqdefa e^{\tfrac{\tau}{2}}\vv(e^{\tau},e^{\tfrac{\tau}{2}}\xxi,z).
\]
Then, from \eqref{eq 1.1}, the pair $(\bH, \bV)$ satisfies
\begin{equation}\label{eq 1.4}
\left\{
\begin{aligned}
&\partial_\tau \bH 
   + \bV \cdot \overline{\nabla} \bH 
   - \bH \cdot \overline{\nabla} \bV
   = (L + e^{\tau}\partial_z^2)\bH, 
   \quad (\tau,\xxi,z) \in \R \times \R^2 \times \T, \\[2mm]
& \bV = \overline{\nabla} \times (-\overline{\Delta})^{-1}\bH, 
   \qquad \overline{\mathrm{div}}\,\bH = 0,
\end{aligned}
\right.
\end{equation}
where the rescaled operators are defined by  
\begin{equation}\label{eq 1.5}
L \eqdefa \Delta_{\xxi} + \frac{\xxi}{2}\cdot \nabla_{\!\xxi} + 1,
\qquad
\overline{\nabla} \eqdefa (\partial_{\xi_1},\,\partial_{\xi_2},\,e^{\tfrac{\tau}{2}}\partial_z),
\qquad
\overline{\Delta} \eqdefa \Delta_{\xxi} + e^{\tau}\partial_z^2,
\qquad
\overline{\mathrm{div}} \eqdefa \overline{\nabla}\!\cdot.
\end{equation}
In the self-similar coordinates, the 3D Oseen vortex filament \eqref{eq 1.2} reduces to the equilibrium $(\alpha\oome^{\G}, \alpha\UU^\G)$ of the rescaled system \eqref{eq 1.4}, where
\begin{equation}
  \label{equilibrium}  \begin{split}
        &\oome^{\G}(\tau,\xxi,z)\eqdefa  e^\tau \eeta^{\G}(e^\tau,e^\frac{\tau}{2}\xxi,z) = (0,0,G(\xxi))^\tr,\\
        &\quad \UU^\G(\tau,\xxi,z) \eqdefa e^\frac{\tau}{2} \boldsymbol{\rm v}^\G(e^\tau,e^\frac{\tau}{2}\xxi,z)= \begin{pmatrix}
            \vv^\G(\xxi)\\0
        \end{pmatrix}.
    \end{split}
\end{equation}
To study the long-time behavior of the fully nonlinear Navier--Stokes equations near the Oseen vortex filament, we perform a perturbation decomposition of the rescaled equations \eqref{eq 1.4} around the equilibrium state \eqref{equilibrium}. Setting $\hh(\tau)= \al \oome^\G + \oome(\tau)$, we deduce from \eqref{eq 1.4} that $\oome$ satisfies
\begin{equation}\label{eq 1.7}
\left\{
    \begin{aligned}
    &\pa_\tau \oome -  (L+e^\tau\pa_z^2) \oome+ \al\big(\UU^{\G}\cdot \bgrad \oome+\UU\cdot \bgrad \oome^{\G}-\oome\cdot \bgrad \UU^{\G}-\oome^{\G}\cdot\bgrad \UU\big)\\
    &\qquad \qquad \qquad \qquad \qquad  =- \UU\cdot \bgrad \oome + \oome\cdot \bgrad \UU, \\
    &\UU \eqdefa \bgrad \times (-\blap)^{-1}\oome, \quad \bdive \, \oome=0.
    \end{aligned}
\right.
\end{equation}
Since our analysis focuses on the asymptotic decay of solutions, it is convenient to consider \eqref{eq 1.7} with initial data prescribed at $\tau = 0$ (i.e., at $t = 1$):
\begin{equation}
    \oome(\tau=0)=\oome(0).
\end{equation}
The flow near a 3D Oseen vortex filament exhibits a delicate interplay between vortex stretching, shearing, and mixing. This interplay constitutes a fundamental mechanism underlying the formation and evolution of coherent structures in high-Reynolds-number fluids. From a physical perspective, a vortex filament can be regarded as a thin tube carrying a highly concentrated vorticity distribution, whose axial velocity profile induces strong azimuthal differential rotation.

Non-axisymmetric perturbations transported by this background field are continuously stretched, wound up, and filamented—closely analogous to the spiral vortex layers observed in studies of straight vortex tubes (see, e.g., \cite{Kaw05}). This stretching–shearing process rapidly generates smaller spatial scales, thereby substantially enhancing viscous dissipation and producing decay rates for transverse disturbances that are much faster than those predicted by pure heat diffusion.

In two-dimensional vortex dynamics, analogous shearing-induced mixing mechanisms have been thoroughly investigated, most prominently for the Lamb–Oseen vortex, where azimuthal differential rotation triggers enhanced dissipation and drives the axisymmetrization of non-axisymmetric modes, as rigorously justified in the mathematical work of \cite{Ga18}.

In the three-dimensional case, however, the presence of an axial variable, the resulting anisotropy, and the additional vortex-stretching term introduce significant analytical difficulties. A general framework for handling vortex-filament initial data in the 3D Navier–Stokes equations is developed in \cite{BGH23}. There, the global stability of the Oseen vortex filament is established in a self-similar anisotropic function space, together with local well-posedness for general Jordan-curve vortex filaments. While this work provides a foundational framework for the analysis of viscous vortex-filament dynamics, it does not furnish quantitative information on the long-time decay rates of non-axisymmetric perturbations. Regarding other studies on the well-posedness or uniqueness criteria for vortex filaments, we refer the interested reader to \cite{GH25, BG21, GS19}.

The purpose of the present work is to provide a precise, quantitative description of the long-time dynamics near the Oseen vortex filament in the high-Reynolds-number regime and to identify the specific role played by the shear--mixing mechanism in the asymptotic behavior. Within an anisotropic self-similar coordinate system adapted to the filament geometry, we demonstrate that the combined effect of strong azimuthal shear and axial diffusion yields an optimal decay rate $t^{-\kappa|\alpha|^{1/2}}$ for all non-axisymmetric perturbations.

At the linear level, this phenomenon manifests itself as a sharp spectral lower bound $\Sigma(\alpha) \sim |\alpha|^{1/2}$ for the nonlocal Oseen operator $L_\perp - \alpha \Lambda_\perp$, and we identify an explicit spectral point precisely at this optimal bound. Combined with the spectral lower bound estimates established in \cite{LWZ}, our analysis completely resolves the conjecture proposed by Gallay and Maekawa in \cite{GM} concerning the asymptotic scaling laws of the spectral bound $\Sigma(\alpha)$ and the pseudospectral bound $\Psi(\alpha)$. This provides a rigorous mathematical explanation of the shear--mixing mechanism near the 3D Oseen vortex filament.

Before presenting our main results, we introduce several function spaces that will be used throughout the paper. 
Let $\boldsymbol{\xi} = r(\cos\theta, \sin\theta)$. For any three-dimensional vector field $\oome$, we write
\begin{equation}\label{eq 1.9}
    \oome = \Omega^r \mathbf{e}_r + \Omega^\theta \mathbf{e}_\theta + \Omega^z \mathbf{e}_z = \Omega^1 \mathbf{e}_1 + \Omega^2 \mathbf{e}_2 + \Omega^z \mathbf{e}_z \eqdefa \oome^\sxxi + \Omega^z \mathbf{e}_z,
\end{equation}
where 
\begin{equation*}
    \begin{pmatrix}
        \mathbf{e}_r \\ \mathbf{e}_\theta 
    \end{pmatrix} = \begin{pmatrix}
        \cos \theta & \sin \theta \\
        -\sin \theta & \cos \theta 
    \end{pmatrix}
   \begin{pmatrix}
        \mathbf{e}_1\\   \mathbf{e}_2 
    \end{pmatrix} 
     \eqdefa \cR_\theta\begin{pmatrix}
        \mathbf{e}_1 \\ \mathbf{e}_2 
    \end{pmatrix}, \quad \text{and therefore} \quad \begin{pmatrix}
        \Omega^r \\ \Omega^\theta 
    \end{pmatrix} = \cR_\theta \begin{pmatrix}
        \Omega^1\\ \Omega^2
    \end{pmatrix}.
\end{equation*}
We adopt the convention that $\hat{f}$ denotes the Fourier transform of $f$ with respect to the $z$-variable only; that is,
\[
\hat{f}(\zeta) \eqdefa \mathcal{F}_{z \to \zeta}(f)(\zeta) \eqdefa (2\pi)^{-\f12} \int_{\T} f(z)e^{-iz\zeta} dz.
\] 
Similarly, the Fourier projection in the $\theta$-variable is defined by
\begin{equation}
\begin{split}
   & \cP^\theta_k f \eqdefa {(2\pi)^{-\f12}}\int_\T f(\theta) e^{-ik\theta} d\theta ,\qquad \cP^\theta_{\neq} f \eqdefa f(\theta)-\cP^\theta_0 f.
\end{split}
\end{equation}
When there is no risk of ambiguity, $\mathcal{P}^\theta_k f$ and $\mathcal{P}^\theta_{\neq} f$ will be abbreviated as $f_k$ and $f_{\neq}$, respectively.

For $\boldsymbol{\oome}_{\neq}$ and $\boldsymbol{\oome}_0$, their definitions are, respectively,
\[
\boldsymbol{\oome}_{\neq} \eqdefa 
\Omega^r_{\neq}\mathbf{e}_r
+ \Omega^\theta_{\neq}\mathbf{e}_\theta
+ \Omega^z_{\neq}\mathbf{e}_z,
\qquad
\boldsymbol{\oome}_0 \eqdefa
\Omega^r_{0}\mathbf{e}_r
+ \Omega^\theta_{0}\mathbf{e}_\theta
+ \Omega^z_{0}\mathbf{e}_z.
\]

We now introduce a collection of weighted function spaces that are fundamental to the subsequent analysis:
\begin{equation}\begin{split}\label{eq 1.13}
L^2&(m) \eqdefa  \left\{ w : \  \|w\|_{L^2(m)}^2 \eqdefa  \scalebox{1.2}{$\int_{\mathbb{R}^2}$} (1 + |\xxi|^2)^m |w(\xxi)|^2 \d\xxi  < \infty \right\}, \\
L^p&(\oo)\eqdefa \left\{w :\   \|w\|_{L^p(\oo)}\eqdefa  \|G^{-1/2}w\|_{L^p(\R^2;d \scalebox{0.8}{$\xxi$})}<+\oo \right\},\\
 Y &\eqdefa L^2(\oo) , \qquad 
Y_0 \eqdefa \left\{ w \in Y: \  \scalebox{1.2}{$\int_{\mathbb{R}^2}$} w(\xxi) \d\xxi = 0 \right\} = \{G\}^{\perp}, \\
Y_1 &\eqdefa \left\{ w \in Y_0:\  \scalebox{1.2}{$\int_{\mathbb{R}^2}$} \xxi w(\xxi)\d\xxi = 0 \right\} = \{G, \partial_1 G, \partial_2 G\}^{\perp}, \\
Y_2 &\eqdefa \left\{ w \in Y_1:\  \scalebox{1.2}{$\int_{\mathbb{R}^2}$} |\xxi|^2 w(\xxi)\d\xxi = 0 \right\} = \{G, \partial_1 G, \partial_2 G, \lapxi G\}^{\perp}.
\end{split}\end{equation}
Here, ``$\perp$'' denotes orthogonality in the space $Y$ with respect to its natural inner product $\langle \cdot, \cdot \rangle_Y$. 
Furthermore, the space $Y$ admits the following orthogonal frequency decomposition:
$$
Y= \ooplus_{n \in \mathbb{N}} X_n  \with X_n\eqdefa \bigl\{w\in Y:\ w(\xxi)= a(r)\cos(n\theta) + b(r)\sin(n\theta) \bigr\}, 
$$
with $X_j \perp X_k$ for all $j \neq k$.

Due to the translational invariance of the Oseen vortex filament in the $z$-direction, the analysis of perturbations naturally employs an anisotropic functional framework. In particular, higher regularity is imposed in the axial $z$-direction than in the horizontal variables. Following \cite{BGH23}, we introduce the $X$-valued Wiener algebra $B_z(X)$, equipped with the norm
\begin{equation}\label{definition Bz}
    \|f\|_{B_z(X)} \eqdefa
        \sum_{\zeta\in\Z} \|\hat{f}(\boldsymbol{\cdot},\zeta)\|_X.
\end{equation}

Finally, for convenience, we define a subspace $Z$ of $B_z(Y)$, which is the space to which the component $\Omega^z$ belongs:
\begin{equation}
    Z\eqdefa \left\{ w \in B_z(Y) :\  \scalebox{1.2}{$\int_{\mathbb{R}^2}$} \scalebox{1.2}{$\int_{\T}$}\ \xxi w(\xxi) \d\xxi dz =0, \quad \scalebox{1.2}{$\int_{\mathbb{R}^2}$} \scalebox{1.2}{$\int_{\T}$}\  w(\xxi)\d\xxi dz =0 \right\}.
\end{equation}
We then have the decomposition
\begin{equation}
Z =  Z_0\ooplus Z_{\neq}  \with   Z_0 \eqdefa Z\cap B_z(X_0) \andf  Z_{\neq} \eqdefa Z\cap B_z(X_{\neq}),
\end{equation}
where $X_{\neq}\eqdefa \bigcup_{k\neq 0 }X_k.$

\subsection{Main theorem}
We now present the main result of this paper concerning the nonlinear asymptotic stability of the 3D Oseen vortex filament.
\begin{Theorem}\label{Thm 1}
{\sl Assume that the initial data $\oome(0)\in B_z(Y)$ satisfy
\begin{equation}\label{eq 1.12}
\bdive\, \oome(0)=0, \quad
\int_{\R^2\times\T} \Omega^z(0) \, d\xxi dz=0, \quad
\int_{\R^2\times\T} \xxi \Omega^z(0) \, d\xxi dz=0.
\end{equation}
There exist universal constants $\kappa,\alpha_0>0$ and  $\epsilon(\alpha), C_\alpha>0$, depending only on $\alpha\geq \alpha_0$, such that if
\begin{equation}
\|\oome(0)\|_{ B_z(Y)} \leq \epsilon(\al),
\end{equation}
then the system \eqref{eq 1.7} admits a unique solution
$\oome(\tau)\in C([0,\oo); B_z(Y))$ satisfying, for all $\tau\ge 0$,
\begin{equation}\label{S1eq20}
\begin{split}
& \bdive\, \oome(\tau)=0, \quad
\int_{\R^2\times\T} \Omega^z(\tau) \, d\xxi dz=0, \quad
\int_{\R^2\times\T} \xxi \Omega^z(\tau) \, d\xxi dz=0,
\\
& (1+\tau)^{-1} e^{\tau/2} \|\mathcal{M}_\tau \oome_{0}(\tau)\|_{B_z(Y)}
+ e^{\kappa|\al|^{1/2}\tau} \|\mathcal{M}_\tau \oome_{\neq}(\tau)\|_{B_z(Y)}
\le C_\al\, \|\oome(0)\|_{B_z(Y)},
\end{split}
\end{equation}
where $\mathcal{M}_\tau \eqdefa \exp\bigl((e^{\tau/2}-1)|D_z|\bigr)$ is a Fourier multiplier in the $z$-variable with symbol $\exp\bigl((e^{\tau/2}-1)|\zeta|\bigr).$}
\end{Theorem}

Translating the above theorem back to the original physical $(t,x)$ coordinates yields the following result.
\begin{Theorem}
  {\sl  Let $\eeta^G$  be given by \eqref{eq 1.2},
 we assume that $e^{|x|^2/8} (\eeta-\alpha \eeta^G)(1)\in {B_z(L^2_\xx)}$ and
 \begin{equation*}
 \dive \eeta(1)=0,\quad \int_{\R^2\times\T}\eta^z(1,x,z)\,dx\,dz=\alpha \andf \int_{\R^2\times\T}x\eta^z(1,x,z)\,dx\,dz= 0.
 \end{equation*}
 Then there exist universal constants $\kappa,\alpha_0>0$ and  $\epsilon(\alpha), C_\alpha>0$, depending only on $\alpha\geq \alpha_0$, such that if
\begin{equation*}
 \big\|e^{|x|^2/8} (\eeta-\alpha \eeta^G)(1)\big\|_{B_z(L^2_\xx)}\leq \epsilon(\al),
\end{equation*} 
 the system  \eqref{eq 1.1}  with data $\eeta(1)$ at time $1$ has a unique solution $\eeta$ for $t\geq 1$, which satisfies
\begin{equation*}
\begin{split}
     (1+\log t)^{-1} t^{\frac{3}{2}-\frac{1}{p}} \|(\eeta-\alpha \eeta^G)_0(t)\|_{B_z(L^p_{\xx})} &+ t^{\kappa|\al|^{1/2}+1-\frac{1}{p}} \|(\eeta-\alpha \eeta^G)_{\neq}(t)\|_{B_z(L^p_{\xx})} \\
     &\leq C_\al \big\|e^{|x|^2/8} (\eeta-\alpha \eeta^G)(1)\big\|_{B_z(L^2_\xx)}, \qquad \forall t\geq 1.
     \end{split}
\end{equation*}}
\end{Theorem}

Let us comment on the vanishing mass and momentum conditions in \eqref{eq 1.12}.
In fact, the last two orthogonality conditions on $\Omega^z(0)$ in \eqref{eq 1.12} can be removed by adjusting the parameter $\alpha$ and using the translation invariance of the 3D Navier--Stokes equations.
To this end, we write the total vorticity field as
\[
\boldsymbol{H}(\tau,\xxi,z)= \al \oome^{\G} + \oome(\tau),
\]
and define
\[
\al_0\eqdefa \int_{\R^2\times\T} \Omega^z(0) \, d\xxi dz,
\qquad
\boldsymbol{\beta}\eqdefa \int_{\R^2\times\T} \xxi \Omega^z(0) \, d\xxi dz.
\]
It is then equivalent to consider the modified field
\[
\boldsymbol{H}_{\rm mod}(\tau,\xxi,z)
\eqdefa
\boldsymbol{H}(\tau,\xxi+\boldsymbol{\beta}_{\rm mod}e^{-\tau/2},z)
\with
\al_{\rm mod}\eqdefa \al+\al_0,
\qquad
\boldsymbol{\beta}_{\rm mod}\eqdefa \boldsymbol{\beta}/\al_{\rm mod},
\]
which also satisfies \eqref{eq 1.4}.
Consequently,
\[
\oome_{\rm mod}\eqdefa
\boldsymbol{H}_{\rm mod}-\al_{\rm mod}\oome^{\G}
\]
satisfies the orthogonality conditions in \eqref{eq 1.12} and solves \eqref{eq 1.7} with $\al$ replaced by $\al_{\rm mod}$.

Recalling \cite{Ga18}, Gallay used the optimal pseudospectral bound obtained in \cite{LWZ} to derive an asymptotic decay rate of $e^{-c|\alpha|^{1/3}\tau}$ for the non-zero-frequency part of the perturbation near the 2D Lamb–Oseen vortex; see Theorem \ref{arma-18} below. This decay rate corresponds to the so-called enhanced dissipation mechanism, which results from the stretching effect induced by the underlying shear flow. Furthermore, Gallay \cite{Ga18} suggested that the eventual asymptotic rate might be improved to $e^{-\kappa|\alpha|^{1/2}\tau}$. The present paper not only provides a rigorous proof of this improved decay rate but also demonstrates that this rate is, in fact, optimal.

We begin by rewriting equation \eqref{eq 1.7} in self-similar cylindrical coordinates:
\begin{subequations}\label{eq 1.14}
\begin{gather}
    \partial_\tau
    \begin{pmatrix}
        \Omega^r \\ \Omega^\theta
    \end{pmatrix}
    - \bigl(e^\tau \partial_z^2 + \sL - \alpha \Pi\bigr)
    \begin{pmatrix}
        \Omega^r \\ \Omega^\theta
    \end{pmatrix}
    - \alpha
    \begin{pmatrix}
        R^r(\oome) \\ R^\theta(\oome)
    \end{pmatrix}
    =
    - \begin{pmatrix}
        B^r[\UU,\oome] \\ B^\theta[\UU,\oome]
    \end{pmatrix}, \\[1mm]
    \partial_\tau \Omega^z
    - \bigl(e^\tau \partial_z^2 + L - \alpha \Lambda\bigr)\Omega^z
    - \alpha R^z(\oome)
    = - B^z[\UU,\oome].
\end{gather}
\end{subequations}
Here, $\boldsymbol{R}(\oome)\eqdefa
R^r \mathbf{e}_r + R^\theta \mathbf{e}_\theta + R^z \mathbf{e}_z$ with
\begin{equation}\label{eq 1.15}
\begin{split}
R^r(\oome) &\eqdefa G e^{\tau/2} \partial_z U^r, \qquad
R^\theta(\oome) \eqdefa G e^{\tau/2} \partial_z U^\theta, \\
R^z(\oome) &\eqdefa
G e^{\tau/2} \partial_z U^z
- \nabla_{\!\xi} G \cdot
\bigl(\UU^\sxxi - \nabla_{\!\xi}^\perp \Delta_{\xi}^{-1}\Omega^z\bigr).
\end{split}
\end{equation}
The bilinear terms are defined by
\[
\boldsymbol{B}(\UU,\oome)\eqdefa
\UU\!\cdot\!\bgrad \oome - \oome\!\cdot\!\bgrad \UU
\eqdefa
B^r \mathbf{e}_r + B^\theta \mathbf{e}_\theta + B^z \mathbf{e}_z,
\]
where
\begin{equation}\label{eq 1.16}
\begin{split}
B^r[\UU,\oome] &= \UU\!\cdot\!\bgrad \Omega^r - \oome\!\cdot\!\bgrad U^r, \qquad
B^z[\UU,\oome] = \UU\!\cdot\!\bgrad \Omega^z - \oome\!\cdot\!\bgrad U^z, \\
B^\theta[\UU,\oome] &=
\UU\!\cdot\!\bgrad \Omega^\theta - \oome\!\cdot\!\bgrad U^\theta
+ r^{-1}(U^\theta \Omega^r - \Omega^\theta U^r).
\end{split}
\end{equation}
The rescaled operators are given by
\begin{equation*}
\begin{split}
&\Delta_{\sxxi} = \partial_r^2 + \frac{1}{r}\partial_r + \frac{1}{r^2}\partial_\theta^2, \qquad
L \eqdefa \Delta_{\sxxi} + \frac{\xxi}{2}\!\cdot\!\nabla_{\!\sxxi} + 1, \qquad\Lambda f \eqdefa
\boldsymbol{v}^{\G}\!\cdot\!\nabla_{\sxxi} f
+ \nabla_{\sxxi}^\perp \Delta_{\sxxi}^{-1} f \cdot \nabla_{\sxxi} G, \\[1mm]
&\sL \eqdefa
\mathcal{R}_\theta \circ
\begin{pmatrix}
L & 0 \\ 0 & L
\end{pmatrix}
\circ \mathcal{R}_\theta^{\!\tr}
=
(L - r^{-2})\mathrm{Id}
- \frac{2}{r^2}
\begin{pmatrix}
0 & 1 \\ -1 & 0
\end{pmatrix}\partial_\theta,
\qquad
\mathrm{Id} \eqdefa
\begin{pmatrix}
1 & 0 \\ 0 & 1
\end{pmatrix}, \\[1mm]
&\Pi \eqdefa
\begin{pmatrix}
S(r)\partial_\theta & 0 \\
- r S'(r) & S(r)\partial_\theta
\end{pmatrix},
\qquad
S(r)\eqdefa \frac{1}{2\pi r^2}\bigl(1 - e^{-r^2/4}\bigr).
\end{split}
\end{equation*}

It is worth noting that the operator $L - \alpha\Lambda$ is commonly referred to as the linearized Oseen vortex operator. Its spectral properties fundamentally govern the long-time behavior of solutions near both the 2D Lamb–Oseen vortex and the 3D Oseen vortex filament. More detailed spectral properties of $L - \alpha\Lambda$ will be discussed in Section~\ref{sect.2}. Importantly, $L$ is self-adjoint in the space $Y$, whereas $\Lambda$ is skew-adjoint in $Y$, and we have
\begin{equation}
{\rm Ker}\,\Lambda = X_0 \ooplus \operatorname{span}\{\partial_1 G,\, \partial_2 G\},
\qquad
({\rm Ker}\,\Lambda)^{\perp} = Y_1 \ccap \ooplus_{n\neq 0} X_n.
\end{equation}

Beyond its crucial role in characterizing the long-time dynamics of the fluid equations, the operator $L - \alpha\Lambda$ also serves as a canonical example in the spectral theory of non-self-adjoint operators. As the parameter $\alpha$ becomes large, it represents a transition from the self-adjoint operator $L$ to the non-self-adjoint operator $L - \alpha\Lambda$. This transition is driven by the skew-adjoint perturbation $\alpha\Lambda$, which is relatively compact with respect to $L$. Since the primary interest lies in understanding how the parameter $\alpha$ influences the non-self-adjoint nature of the operator, it is natural to investigate the properties of $L - \alpha\Lambda$ restricted to the space $({\rm Ker}\,\Lambda)^{\perp}$; namely,
\[
L_\perp - \alpha \Lambda_\perp \eqdefa (L - \alpha\Lambda)\big|_{({\rm Ker}\,\Lambda)^{\perp}}.
\]

To accurately quantify the influence of the rapid rotation parameter $\alpha$, it is customary to introduce the following two quantities, which characterize the spectral properties:
\begin{equation}\label{spectral bound def}
    \begin{split}
& \textit{spectral lower bound}: \quad \Sigma(\alpha) \eqdefa \inf\Bigl\{\Re(z)\, : \, z \in
  \sigma(-L_\perp + \alpha\Lambda_\perp)\Bigr\}, \\
& \textit{pseudospectral bound}: \quad \Psi(\alpha) \eqdefa  \Bigl(\sup_{\la \in \mathbb{R}}\|(L_\perp
  -\alpha\Lambda_\perp-i\la)^{-1}\|_{Y \to Y}\Bigr)^{-1}.
    \end{split}
\end{equation}
Accordingly, Gallay and Maekawa proposed the following interesting conjecture in \cite{GM}:
\begin{equation}\label{eq conjecture}
\textbf{Conjecture: } \Sigma(\alpha) =
\mathcal{O}(|\alpha|^{1/2}) \text{ and } \Psi(\alpha) = \mathcal{O}(|\alpha|^{1/3}) \text{ as }
|\alpha| \to \infty.
\end{equation}
In earlier work \cite{LWZ}, Li, Wei, and Zhang established the following result.

\begin{Theorem}[\cite{LWZ}]\label{Thm 1.3}
{\sl There exists a constant $C>0$, independent of $\alpha$, such that as $|\alpha|\to\infty$,
\begin{equation}
    C^{-1}|\alpha|^{1/2} \leq \Sigma(\alpha),
    \qquad
    C^{-1}|\alpha|^{1/3} \le \Psi(\alpha) \le C|\alpha|^{1/3}.
\end{equation}}
\end{Theorem}
A more detailed review of results related to Conjecture~\eqref{eq conjecture}
is provided in Section~\ref{sect.2}.
The following theorem states the second main result of the paper,
which fully confirms the first part of Conjecture~\eqref{eq conjecture}.
Indeed, by combining Proposition~\ref{Prop 5.3}, Proposition~\ref{Prop 5.6}, and Theorem~\ref{Thm 1.3}, we obtain
\begin{Theorem}\label{Thm 2}
{\sl There exist constants $C_1, C_2 > 0$, independent of $\alpha$, such that as $|\alpha|\to+\infty$,
\begin{equation}
    C_1 |\alpha|^{1/2} \le \Sigma(\alpha) \le C_2 |\alpha|^{1/2}.
\end{equation}}
\end{Theorem}

The proof of Theorem~\ref{Thm 2} draws inspiration from the works \cite{GGN,HS96}. The basic strategy is to construct a family of complexified harmonic oscillator operators to approximate
$L_k - \alpha\Lambda_k = (L_\perp - \alpha \Lambda_\perp)\big|_{X_k}$.
By carefully implementing the complex deformation method, we identify a sequence of spectral points of $L_k - \alpha \Lambda_k$ whose real parts lie at the scale $-\mathcal{O}(|\alpha|^{1/2})$, thereby establishing Theorem~\ref{Thm 2}. It is particularly worth emphasizing that this approach can directly handle operators with a nonlocal structure—the presence of nonlocal terms in non-self-adjoint operators constitutes one of the most ubiquitous and fundamentally challenging difficulties in the spectral analysis of fluid stability theory.

\noindent{\bf Remarks on the optimality of $e^{-\kappa|\alpha|^{1/2}\tau}$:}
Here we explain how the optimality of the spectral bound in Theorem~\ref{Thm 2} leads to the optimality of the asymptotic decay rate in Theorem~\ref{Thm 1}. To this end, consider an initial perturbation $\oome(0)$ satisfying
\begin{equation}\label{eq 1.31}
\Omega^r(0,\xxi,z)=\Omega^\theta(0,\xxi,z)=0 \quad \text{and} \quad \Omega^z(0,\xxi,z)=\Omega^z(0,\xxi) \in Y_1.
\end{equation}
Then \eqref{eq 1.31} holds for all $\tau \ge 0$, and the 3D Navier--Stokes equations reduce to a 2D Navier--Stokes system for $\Omega^z$. In particular, $\Omega^z_{\neq}(\tau) \in ({\rm Ker}\,\Lambda)^{\perp}$, and the linear part of the equation for $\Omega^z_{\neq}$ becomes
\begin{equation*}
\partial_\tau \Omega^z_{\neq} - (L_{\perp}-\alpha \Lambda_{\perp}) \Omega^z_{\neq} = 0.
\end{equation*}
By Theorem~\ref{Thm 2} and relatively compact perturbation theory, we obtain
\begin{equation*}
\lim_{\tau\to +\infty} \tau^{-1} \ln\big\|e^{\tau(L_{\perp}-\alpha \Lambda_{\perp})}\big\|_{({\rm Ker}\,\Lambda)^{\perp}}
= -\Sigma(\alpha) \ge -C_2 |\alpha|^{1/2},
\end{equation*}
which shows that the decay rate $e^{-\kappa |\alpha|^{1/2} \tau}$ is optimal at the linear level. As for the nonlinear evolution, when the initial perturbation is sufficiently small, the nonlinear terms do not dominate, and hence the rate $e^{-\kappa|\alpha|^{1/2}\tau}$ remains optimal in the nonlinear regime as well. We will not pursue this further here.

\subsection{Notations}\label{Notations}
Throughout the paper, we fix the coefficient
\[
a(\tau) = 1 - e^{-\tau}.
\]
Unless otherwise specified, $C$ denotes a universal constant, while $C_\alpha$ denotes a constant depending on the parameter $\alpha$. We write $g \lesssim f$ if there exists a constant $C>0$ such that $g \le Cf$; when it is important to indicate that the implicit constant depends on $\alpha$, we write $g \lesssim_\alpha f$. Similarly, we write $g \thicksim f$ if both $g \lesssim f$ and $f \lesssim g$ hold.

Let $\ltr \eqdefa L^2(\R^+;dr)$ denote the Hilbert space equipped with its natural inner product $\langle \cdot, \cdot \rangle_{\ltr}$. For any unbounded, densely defined, closed operator $A$ on a Banach space $X$ with domain $D(A)$, we denote by $\sigma_X(A)$ the spectrum of $A$ and define its spectral bound as
\[
s_X(A) \eqdefa \sup_{z \in \sigma_X(A)} \mathrm{Re}\, z.
\]
For any bounded linear operator $A \in \mathcal{L}(X)$, we abbreviate the operator norm $\|A\|_{X \to X}$ simply as $\|A\|_X$.

\renewcommand{\theequation}{\thesection.\arabic{equation}}
\setcounter{equation}{0}
\section{Asymptotic stability results for the 2D Lamb--Oseen vortex}\label{sect.2}
In this section, we present a brief survey of stability problems near the Lamb--Oseen vortex, whose dynamics are governed by the two-dimensional Navier--Stokes equations:
\begin{equation}\label{eq 2.1}
\partial_t \omega - \Delta \omega + \uu \cdot \nabla \omega = 0,
\qquad (t,\boldsymbol{x}) \in \R^+ \times \mathbb{R}^2,
\end{equation}
where the vorticity $\omega$ is a scalar function and the associated velocity field $\uu(\xx)$ is determined by the Biot--Savart law:
\begin{equation}\label{eq 2.2}
\uu(t,\boldsymbol{x})
= \frac{1}{2\pi} \int_{\mathbb{R}^2}
\frac{(\boldsymbol{x} - \boldsymbol{y})^{\perp}}
{|\boldsymbol{x} - \boldsymbol{y}|^2}
\, \omega(t,\boldsymbol{y})\, d\boldsymbol{y}.
\end{equation}

In particular, when the initial vorticity satisfies $\omega(t=0)=\alpha\delta_0$, the exact solution is the so-called Lamb--Oseen vortex:
\begin{equation}\label{eq 2.3}
\omega(t,\boldsymbol{x})
= \frac{\alpha}{t}\,
G\!\left(\frac{\boldsymbol{x}}{\sqrt{t}}\right),
\qquad
\uu(t,\boldsymbol{x})
= \frac{\alpha}{\sqrt{t}}\,
\vv^{\G}\!\left(\frac{\xx}{\sqrt{t}}\right),
\end{equation}
which possesses a self-similar structure. Here, the circulation number $\alpha \gg 1$ measures the strength of rotation.

The well-posedness of \eqref{eq 2.1} with initial data in $\mathcal{M}(\R^2)$ (the space of finite signed measures) has been studied in a series of works
\cite{ben-artzi:1994,brezis:1994,kato:1994,GGL05,giga:1988,GG05}.

\begin{Theorem}[\cite{GG05}]
{\sl For any $\mu \in \mathcal{M}(\R^2)$, equation \eqref{eq 2.1}
with initial data $\mu$ has a unique solution $\omega(t)\in C\big((0,
\infty),\, L^1_\xx(\R^2)\cap L^\infty_\xx(\R^2)\big)$
satisfying
\begin{equation*}
    \sup_{t>0}\|\omega(t)\|_{L^1_\xx(\R^2)} \le M < \infty,
    \qquad
    \omega(t) \rightharpoonup \mu \ \text{as} \ t \to 0.
\end{equation*}}
\end{Theorem}

A particularly remarkable feature of the Lamb--Oseen vortex is that any solution of \eqref{eq 2.1} with measure-valued initial data (and, thanks to the smoothing effect of \eqref{eq 2.1}, it suffices to consider initial data in $L^1(\R^2)$) eventually converges to a Lamb--Oseen vortex of some circulation strength; see \cite{giga:1988b,GW2}.
\begin{Theorem}[\cite{GW2}]\label{Thm 2.2}
  {\sl For any $\omega_0\in L^1_\xx(\R^2)$, the solution $\omega(t,\xx)$ of \eqref{eq 2.1} satisfies
\begin{equation*}
\begin{split}
    \lim_{t\rightarrow +\oo}  t^{1-\frac{1}{p}}\Bigl\|\omega(t,\cdot)- \frac{\al}{ t}G\Bigl(\frac{\cdot}{\sqrt{ t}}\Bigr) \Bigr\|_{L^p_\xx} =0 \quad  \text{for } \ 1\leq p\leq +\oo,\\
     \lim_{t \to \infty}  t^{\f12 - \frac{1}{q}} \Big\|
  \uu(t,\cdot) - \frac{\al}{\sqrt{ t}}\vv^{\G}\Bigl(\frac{\cdot}{\sqrt{ t}}\Bigr)
  \Big\|_{L^q_\xx} = 0 \quad
  \text{ for}\quad2 < q \le \infty,
\end{split}
\end{equation*}
where $\al\eqdefa \int_{\R^2} w_0(\xx) \d\xx$.}
\end{Theorem}

In contrast to the three-dimensional case---where energy injected into the flow at large scales is transferred to progressively smaller scales until it is dissipated by viscosity---experimental and numerical studies in the two-dimensional case indicate that even in highly turbulent, high-Reynolds-number flows, smaller vortices tend to merge and form increasingly large coherent structures. Theorem~\ref{Thm 2.2} shows that, in the full space $\R^2$, this process continues until only a single vortex remains, in the $L^1_\xx$ sense. Moreover, as presented in the following subsections, Gallay and Wayne \cite{GW2} initiated the asymptotic analysis of \eqref{eq 2.1} in stronger topologies---namely, in the spaces $L^2(m)$ and $L^2(\infty)$. In these settings, the Lamb--Oseen vortices are equilibria, and the vorticity profiles enjoy significantly stronger spatial localization than merely belonging to $L^1$.

Recall the self-similar change of variables,
\[
\xxi = \xx/\sqrt{t}, \qquad \tau = \log t,
\]
and define the self-similar profiles of $w$ and $\boldsymbol{v}$ by
\begin{equation}\label{eq 2.4}
w(\tau,\xxi) = t\, \omega(t,\xx),
\qquad
\boldsymbol{v}(\tau,\xxi) = t^{1/2}\, \uu(t,\xx) .
\end{equation}
Then we obtain
\begin{equation}\label{eq 2.5}
\partial_\tau w + \vv \cdot \gradxi w =\bigl(\Delta_{\boldsymbol{\xi}} + \frac{\boldsymbol{\xi}}{2}\cdot \nabla_{\!\boldsymbol{\xi}} + 1\bigr)w=
 L w,
\qquad
\vv \eqdefa \gradxi^{\perp}\lapxi^{-1} w.
\end{equation}
The well-posedness of \eqref{eq 2.5} was established in \cite[Theorem~3.2]{gallay:2001a}.
By construction, the family of Lamb--Oseen vortices $\{\alpha G\}_{\alpha\in\R}$ consists of equilibria of \eqref{eq 2.5}.
Perturbations near $\alpha G$ satisfy
\begin{equation}\label{NS2Dwper}
\partial_\tau \tilde w + \tilde{\vv} \cdot \nabla \tilde w
= (L - \alpha \Lambda)\tilde w ,
\end{equation}
where the nonlocal linear operator $\Lambda$ is defined by
\begin{equation}\label{Lamdef}
\Lambda \tilde w \eqdefa
\vv^{\G} \cdot \gradxi \tilde w
+ \tilde{\vv} \cdot \gradxi G,
\qquad
\tilde{\vv} = \gradxi^{\perp}\lapxi^{-1} \tilde w .
\end{equation}
\subsection{Asymptotic stability in $L^2(m)$ for $ m>1$}\label{sect. 2.1}
Note that, by the definition of $L^2(m)$ in \eqref{eq 1.13} and Hölder's inequality, we have the continuous embedding
$L^2(m)\hookrightarrow L^1_\sxxi(\R^2)$.
In this section, the operators $L$ and $\La$ are regarded as closed operators on $L^2(m)$ with domains
\begin{equation*}
\begin{split}
    D_m(L) &\eqdefa \big\{ w\in L^2(m) \,:\, \lapxi w,\ \xxi\!\cdot\!\gradxi w \in L^2(m) \big\}, \\
    D_m(\La) &\eqdefa \big\{ w\in L^2(m) \,:\, \vv^{\G}\!\cdot\!\gradxi w \in L^2(m) \big\}.
\end{split}
\end{equation*}

Two invariant subspaces of $L^2(m)$ associated with $L$ and $\La$ are introduced as follows:
\begin{equation}
\begin{split}
    L^2_0(m) &\eqdefa \left\{ w\in L^2(m) \,:\,
    \scalebox{1.2}{$\int_{\mathbb{R}^2}$} w(\xxi)\, d\xxi = 0 \right\},\\
    L^2_1(m) &\eqdefa \Bigl\{ w\in L^2_0(m) \,:\,
    \scalebox{1.2}{$\int_{\mathbb{R}^2}$} \xxi\, w(\xxi)\, d\xxi = 0 \Bigr\}.
\end{split}
\end{equation}
A direct computation shows that
\begin{equation}
\begin{split}
   &\La G = \La \pa_j G = L G = 0, \qquad
   L\pa_j G = -\tfrac12\, \pa_j G, \quad j=1,2,\\
   & L^2(m)
   = L^2_0(m) \ooplus \operatorname{span}\{G\}
   = L^2_1(m) \ooplus \operatorname{span}\{G,\pa_1 G,\pa_2 G\}.
\end{split}
\end{equation}

Below we summarize several basic spectral properties of $L$ and $L - \alpha \La$.

\begin{Theorem}[\cite{gallay:2001a, GW2}]\label{Lem 2.3}
{\sl For any $m \ge 0$ and $\alpha \in \R$, we have
\begin{equation}
    \sigma_{L^2(m)}(L)
    \,=\,
    \Bigl\{ z \in \mathbb{C} \ \Big|\ \Re z \le \frac{1-m}{2} \Bigr\}
    \,\cup\,
    \Bigl\{ -\frac{k}{2} \ \Big|\ k \in \N \Bigr\},
\end{equation}
and
\begin{equation}
    \sigma_A(L - \alpha \Lambda) \subset
    \begin{cases}
        \bigl\{ z \in \mathbb{C} : \Re z \le \max\bigl(-\frac12, \frac{1-m}{2}\bigr) \bigr\},
        & \text{if $m>1$ and $A = L^2_0(m)$},\\[4pt]
        \bigl\{ z \in \mathbb{C}: \Re z \le \max\bigl(-1, \frac{1-m}{2}\bigr) \bigr\},
        & \text{if $m>2$ and $A = L^2_1(m)$}.
    \end{cases}
\end{equation}
Moreover, there exists a constant $C>0$, possibly depending on $\alpha$, such that
\begin{equation}
  \big\| e^{\tau(L - \alpha\Lambda)} \big\|_{A}
  \,\le\,
  \begin{cases}
       C, & \text{if $m>1$ and $A = L^2(m)$},\\[4pt]
       C e^{-\tau/2}, & \text{if $m>2$ and $A = L^2_0(m)$},\\[4pt]
       C e^{-\tau}, & \text{if $m>3$ and $A = L^2_1(m)$}.
  \end{cases}
\end{equation}}
\end{Theorem}

Using the spectral properties established above, Gallay and Wayne proved in \cite{GW2} that the family $\{\alpha G\}$ is asymptotically stable in $L^2(m)$.
\begin{Theorem}[\cite{GW2}]
{\sl Let $w_0 \in L^2(m)$, and let $w(\tau) \in C^0([0,+\infty), L^2(m))$ be the solution of \eqref{eq 2.5} with initial data $w_0$. Then, as $\tau \to +\infty$,
\begin{equation}
 \|w(\tau) - \alpha G\|_{L^2(m)} =
    \begin{cases}
        \mathcal{O}(1), & \text{if $m>1$},\\[4pt]
        \mathcal{O}(e^{-\tau/2}), & \text{if $m>2$},\\[4pt]
        \mathcal{O}(e^{-\tau}), & \text{if $m>3$ and $w_0 - \alpha G \in L^2_1(m)$},
    \end{cases}
\end{equation}
where $\alpha \eqdefa \displaystyle\int_{\mathbb{R}^2} w_0(\xxi)\, d\xxi$.}
\end{Theorem}

Thus, the above result shows that any solution of the vorticity equation that is sufficiently
localized in $L^2(m)\hookrightarrow L^1_\sxxi$ will converge, as
time goes to infinity, toward the Lamb--Oseen vortex,
regardless of the size of $\alpha$. Returning to the original variables,
Theorem~\ref{Thm 2.2} then follows as a corollary.

It is worth noting that, in contrast to this result, the semi-flow defined by the
Navier--Stokes equations does not preserve the spatial localization of
the velocity field \cite{gallay:2001a}. For instance, even if the initial
velocity $\uu_0$ belongs to $L^2(m)$ for some $m>2$, the corresponding
solution $\uu(\cdot,t)$ may in general fail to lie in $L^2(m)$ for $t>0$.

\subsection{Asymptotic stability in $Y=L^2(\infty)$}
In this subsection, the operators $L$ and $\La$ are regarded as closed operators on $Y$, with domains
\begin{equation*}
\begin{split}
D(L) \eqdefa \big\{ w \in Y \,:\, \lapxi w,\ \xxi \cdot \gradxi w \in Y \big\} \andf
D(\La) \eqdefa \big\{ w \in Y \,:\, \vv^{\G} \cdot \gradxi w \in Y \big\}.
\end{split}
\end{equation*}
Similar to Subsection~\ref{sect. 2.1}, we first recall some basic spectral properties of $L$ and $\La$ in $Y$.

\begin{Theorem}[\cite{gallay:2001a,GM, M3}]\label{Prop 2.5}
{\sl The following statements hold:
\begin{enumerate}
    \item[(1)]
    The operator $L$ is self-adjoint on $Y$ with compact resolvent and purely discrete spectrum:
    \begin{equation*}
        \sigma_Y(L) = \Big\{-\frac{n}{2} \ : \ n \in \N \Big\}.
    \end{equation*}
    Moreover, ${\rm Ker}\, L = \operatorname{span}\{G\}$ and
    ${\rm Ker}\, \bigl(\frac12 + L\bigr) = \operatorname{span}\{\pa_1 G, \pa_2 G\}$.
    The eigenvalue $-\frac{n}{2}$ has multiplicity $n+1$, and the corresponding
    eigenfunctions are the Hermite functions of degree $n$.

    \item[(2)]
    The operator $\La$ is skew-adjoint on $Y$, and
    \begin{equation*}
        {\rm Ker}\,\La = X_0 \ooplus \operatorname{span}\big\{\pa_1 G, \pa_2 G\big\},
        \qquad
        ({\rm Ker}\,\La)^\perp  = Y_1 \ccap \ooplus_{n\neq 0} X_n.
    \end{equation*}
\end{enumerate}}
\end{Theorem}

We now consider the spectrum of the operator $L - \alpha \Lambda$ in $Y$,
a problem that is both interesting and challenging.
On one hand, since $\Lambda$ is a relatively compact perturbation of $L$ in $Y$,
classical perturbation theory (see \cite{Kato}) implies that
$L - \alpha \Lambda$ has a compact resolvent in $Y$.
Consequently, the spectrum $\sigma_Y(L - \alpha \Lambda)$
consists of discrete eigenvalues $\{\lambda_n(\alpha)\}_{n \in \mathbb{N}}$
with finite multiplicity.

On the other hand, although $\Lambda$ is relatively compact with respect to $L$,
the behavior of $L - \alpha \Lambda$ undergoes a transition as $\alpha$ increases from zero,
from that of a self-adjoint operator to that of a skew-adjoint operator.
This suggests that the eigenvalues $\lambda_n(\alpha)$ may gradually tilt toward the imaginary axis
as $|\alpha| \to \infty$.

We first recall the following result.

\begin{Theorem}[\cite{GM,LWZ}]
{\sl For any $\alpha \in \mathbb{R}$, the following spectral inclusions hold:
\begin{equation*}
\sigma_Y(L-\alpha\Lambda)\subset \big\{z\in\mathbb{C} : \Re z \le 0 \big\},
\qquad
\sigma_{Y_0}(L-\alpha\Lambda)\subset \big\{z\in\mathbb{C} : \Re z \le -1/2 \big\},
\end{equation*}
\begin{equation*}
\sigma_{Y_1}(L-\alpha\Lambda)\subset \big\{z\in\mathbb{C} : \Re z \le -1 \big\},
\qquad
\sigma_{Y_2}(L-\alpha\Lambda)\subset \big\{z\in\mathbb{C} : \Re z < -1 \big\}.
\end{equation*}}
\end{Theorem}

Based on these spectral properties, the authors established certain asymptotic behaviors in $L^2(\infty)$ (see \cite{Ga2,GR,GW2}).

\begin{Theorem}[\cite{Ga2,GR,GW2}]
{\sl Let $w_0 \in Y$. Then there exists a unique solution $w(\tau) \in C^0([0,\infty);Y)$ of \eqref{eq 2.5} with initial data $w(0) = w_0$, which satisfies
\begin{equation*}
    \|w(\tau) - \alpha G\|_{Y} \to 0 \quad \text{as } \tau \to +\infty,
\end{equation*}
where $\alpha \eqdefa \displaystyle\int_{\mathbb{R}^2} w_0(\xxi)\,\d \xxi$. Moreover, there exists $\varepsilon > 0$ such that for all $\alpha \in \mathbb{R}$, if
\begin{equation*}
    \|w_0 - \alpha G\|_{Y} \leq \varepsilon,
\end{equation*}
then
\begin{equation}\label{wlocstab}
  \|w(\tau) - \alpha G \|_{Y} \le \min\bigl(1,\,2e^{-\tau/2}\bigr)
  \|w_0 - \alpha G \|_{Y}, \qquad \forall \tau \ge 0.
\end{equation}}
\end{Theorem}

Moreover, for any $\alpha \in \mathbb{R}$, on the subspace ${\rm Ker}\,\Lambda$, the spectrum of the operator $L - \alpha \Lambda$ (which reduces to $L$) is independent of the circulation parameter $\alpha$. In contrast, on $({\rm Ker}\,\Lambda)^\perp$—an invariant subspace of $L - \alpha \Lambda$—the strong shear effect generated by $\alpha \Lambda$ for large $\alpha$ causes the semigroup $e^{-\tau(L-\alpha\Lambda)}$ to exhibit a decay rate stronger than that of pure heat dissipation. This motivates the introduction of the spectral lower bound and pseudospectral bound defined in \eqref{spectral bound def}.

In the definition of $\Sigma(\alpha)$, the sign of the linearized operator is reversed to obtain a positive quantity. Consequently, one immediately has
\[
\Sigma(\alpha) \ge \Psi(\alpha) \ge 1.
\]

In the high-Reynolds-number regime, this stabilizing effect on $({\rm Ker}\,\Lambda)^\perp$ is qualitatively illustrated by the following result:

\begin{Theorem}[\cite{M3}]\label{stabeff}
{\sl One has $\Psi(\alpha) \to \infty$ and $\Sigma(\alpha) \to \infty$
as $|\alpha| \to \infty$. }
\end{Theorem}

In fact, there are good reasons to conjecture, as proposed in \cite{GM}, that:
\begin{equation*}
\textbf{Conjecture: } \Sigma(\alpha) =
\mathcal{O}(|\alpha|^{1/2}) \quad \text{and} \quad \Psi(\alpha) = \mathcal{O}(|\alpha|^{1/3}) \quad \text{as} \quad
|\alpha| \to \infty.
\end{equation*}
First, extensive numerical computations by Prochazka and Pullin \cite{PP1,PP2} indicate that, as $|\alpha| \to \infty$, 
$\Sigma(\alpha) = \mathcal{O}(|\alpha|^{1/2})$.
In addition, rigorous mathematical analysis of a toy model in \cite{GGN} further supports this conjecture. In particular, for the simplified linearized operator $L - \alpha \tilde{\Lambda}$ (i.e., the version where the nonlocal part of $\Lambda$ is omitted), Deng proved in \cite{De1} that, as $|\alpha| \to \infty$, 
$\Psi(\alpha) = \mathcal{O}(|\alpha|^{1/3})$.
The same result also holds for the full linearized operator $L - \alpha \Lambda$ restricted to a subspace smaller than $({\rm Ker}\,\Lambda)^\perp$, where a finite number of Fourier modes in the angular variable (in polar coordinates) are removed; see \cite{De2} for details. In \cite{LWZ}, most of this conjecture was resolved (see Theorem~\ref{Thm 1.3}).

Based on the results in \cite{LWZ}, the nonlinear stability analysis carried out in \cite{Ga18} shows that, in the high-Reynolds-number limit $|\alpha| \to \infty$, the basin of attraction of the Lamb--Oseen vortex becomes very large, and perturbations relax to axisymmetry on a much shorter timescale.

\begin{Theorem}[\cite{Ga18}]\label{arma-18}
  {\sl There exist positive constants 
$C_1$, $C_2$, and $\kappa$ such that, for all $\alpha \in \mathbb{R}$ and all 
initial data $w_0 \in \alpha G + Y_1$ satisfying 
\begin{equation}
  \|w_0 - \alpha G \|_{L^2(\infty)} \,\le\, \frac{C_1\,(1+|\alpha|)^{1/6}}{
  \log(2+|\alpha|)},
\end{equation}
the unique solution of \eqref{eq 2.5} in $L^2(\infty)$ satisfies, for all $\tau \ge 0$, 
\begin{equation}
  \|w(\tau) - \alpha G \|_{L^2(\infty)} \le C_2e^{-\tau} \|w_0 - 
  \alpha G \|_{L^2(\infty)}, 
\end{equation}
\begin{equation}
  \|(1-P^\theta_0)(w(\tau) - \alpha G) \|_{L^2(\infty)} \le C_2\|w_0 - 
  \alpha G \|_{L^2(\infty)} \exp\Bigl(-\frac{\kappa (1+|\alpha|)^{1/3}\tau}{
  \log(2+|\alpha|)}\Bigr),\label{37}
\end{equation}
where $P^\theta_0$ is the orthogonal projection in $L^2(\infty)$ onto the subspace $X_0$
of all radially symmetric functions. }
\end{Theorem}

\section{Sharp spectral lower bound for $L-\alpha\Lambda$ on $L^2(\infty)$}\label{sect. 5}
This section shows that, for any subspace $X_k$, the Oseen vortex operator $L - \alpha\Lambda$ possesses an eigenvalue whose real part lies in a region of order $|\alpha|^{1/2}$. Before demonstrating this, we begin by recalling the relevant operator notations:
\begin{equation*}
L= \lapxi +  \frac{\xxi}{2} \cdot \gradxi +1,\quad  \Lambda= \vv^{\G} \cdot\gradxi + \gradxi G\cdot\gradxi^\perp \lapxi^{-1},
\end{equation*} 
where 
\begin{equation*}
  G(\xxi) = \frac{1}{4\pi}e^{-|\scalebox{0.8}{\xxi}|^2/4},\quad
  \vv^{\G}(\xxi)=\frac{1}{2\pi}\frac{\xxi^\perp}{|\xxi|^2}
  \Bigl(1 -  e^{-|\scalebox{0.8}{\xxi}|^2/4}\Bigr).
\end{equation*}

We recall that on ${\rm Ker }\Lambda$, the spectrum of $L - \alpha\Lambda \equiv L$ is independent of the circulation parameter~$\alpha$. Thus, the main difficulty lies in studying the operator $L_\perp - \alpha \Lambda_\perp,$ which is defined as the restriction of $L - \alpha\Lambda$ to its orthogonal complement $({\rm Ker}\Lambda)^\perp = Y_1 \cap \bigoplus_{k>0} X_k$ (see Proposition~\ref{Prop 2.5}). Previous results show that the spectrum of $L_\perp - \alpha \Lambda_\perp$ indeed depends on the parameter~$\alpha$ in a nontrivial way; see \cite{LWZ,M3,PP1,PP2}.

By restricting \(L - \alpha \Lambda\) to each subspace \(X_k\), which is naturally isometric to \(L^2(\mathbb{R}^+; G^{-1} r\,dr)\), one obtains:
\begin{equation}\label{l xk}
L_k-\al\Lambda_k \eqdefa \left(L-\alpha\Lambda\right)|_{X_k} = \Big(\partial_r^2 + \frac{1}{r}\partial_r+\frac{r}{2}\partial_r - \frac{k^2}{r^2}+1\Big) - i\al k\Big(S(r)- \f{\partial_r G}{r} \Delta_k^{-1}\Big),
\end{equation}
where $S(r)=\frac{1-e^{-r^2/4}}{2\pi r^2} $, $G(r)=\frac{1}{4\pi}e^{-r^2/4}$, and $\Delta_k = \partial_r^2+\frac{1}{r}\partial_r-\frac{k^2}{r^2}$.  

To make the problem more succinct,  we introduce 
\begin{equation}\label{S5eq1}
\beta_k \eqdefa \frac{k\alpha}{8\pi}, \quad 
\sigma(r) \eqdefa \frac{1-e^{-r^2/4}}{r^2/4}, \quad 
g(r) \eqdefa e^{-r^2/8},
\end{equation}
and define the integral operator
$$
\mathcal{K}_k[f](r) \eqdefa  \int_0^\infty K_k(r,s)f(s)\,ds, \with K_k(r,s)\eqdefa \frac{1}{2|k|} \min\Bigl\{\frac{r}{s},\frac{s}{r}\Bigr\}^{|k|}(rs)^{1/2}.
$$
We then consider an operator on \(\ltr \eqdefa L^2(\mathbb{R}^+;dr)\) that is 
equivalent to \(L_k - \alpha \Lambda_k\):
\[
H_k \eqdefa r^{1/2} g^{-1} (L_k - \alpha \Lambda_k)\, r^{-1/2} g
= \Bigl(\partial_r^2 - \frac{k^2 - \tfrac14}{r^2} - \frac{r^2}{16} + \frac12\Bigr)
- i \beta_k \Bigl(\sigma(r) - g \mathcal{K}_k[g \cdot]\Bigr),
\]
with domains specified by
\begin{equation*}
\begin{aligned}
   D &\eqdefa D(H_k)
   = \Bigl\{ w \in \ltr :\; \partial_r^2 w,\ r^{-2} w,\ r^2 w \in \ltr \Bigr\},
   && |k| \ge 2, \\
   D_1 &\eqdefa D(H_k)
   = \Bigl\{ w \in \ltr :\; r^{1/2}\partial_r^2(r^{-1/2}w),\ 
       r^{1/2}\partial_r(r^{-3/2}w),\ r^2 w \in \ltr \Bigr\},
   && |k| = 1.
\end{aligned}
\end{equation*}
Since \(\overline{H_k} = H_{-k}\), the analysis may be restricted to the case \(k \ge 1\).

\subsection{Spectral lower bound of $L_1 - \al \Lambda_1$}\label{sect. 5.1}
This case is equivalent to studying the spectrum of \(H_{1}\). The main difficulty in handling \(H_{k}\) arises from the nonlocal operator \(\mathcal{K}_k\). However, when \(k = 1\), this nonlocal term can be removed by means of the wave operator introduced in \cite{LWZ}. We define:
\[
Tw(r) \eqdefa w(r) + \frac{g(r)}{\sigma'(r)r^{3/2}} \int_0^r s^{3/2} g(s) w(s)\, ds,
\]
\[
T^t w(r) \eqdefa w(r) + r^{3/2} g(r)\int_r^{+\infty} \frac{g(s) w(s)}{\sigma'(s)s^{3/2}} \, ds,
\]
and
\[
\mathcal{L}_1 u \eqdefa T H_1 T^t u 
= \Bigl(\partial_r^2 - \frac{3}{4r^2} -\frac{r^2}{16} - f(r) + \frac12\Bigr)u 
- i\beta_1 \sigma(r)u ,
\]
where
\begin{equation}\label{S5eq2}
f(r) \eqdefa \frac{2g^4}{(\sigma'(r))^2} + \frac{g^2}{\sigma'(r)}\Bigl(\frac{6}{r}-r\Bigr).
\end{equation}

In what follows, the wave operator \(T\) is restricted to the subspace 
\begin{equation}
    \begin{split}
    \mathcal{V} &\eqdefa \Bigl\{\, w \in \ltr \; : \; \int_0^\infty w(r)\, r^{3/2} g(r)\, dr = 0 \Bigr\}\\
    \end{split}
\end{equation}
and is still denoted by $T$.

\begin{Lemma}[\cite{LWZ}, Lemma 5.2]
It holds that $TT^t = \mathrm{Id}_{L^2_r}$ and $T^t T = \mathrm{Id}_{\mathcal{V}}$.
\end{Lemma}

The above lemma shows that \(T\) is a unitary operator from \(\mathcal{V}\) to \(L^2_r\).
Consequently, we obtain the following important spectral equivalence:
\begin{equation}\label{eq 5.3}
    \sigma_{L^2_r}(\mathcal{L}_1)
    = \sigma_{\mathcal{V}}(H_1)
    = \sigma_{X_1 \cap Y_1}(L_1 - \alpha \Lambda_1).
\end{equation}
At this point, the operator \(\mathcal{L}_1\) contains no nonlocal term, and is therefore 
relatively not difficult to analyze. 
The main result of this subsection is stated below, and its proof will be given at the end of the subsection.

\begin{Proposition}\label{Prop 5.3}
 There exists a constant $C_1 > 0$ such that
\[
s_{X_1\cap Y_1}(L_1 - \al \Lambda_1) = s_{L^2_r}(\mathcal{L}_1) \geq - C_1 |\al|^{1/2}.
\]
\end{Proposition}

The main idea originates from \cite{GGN}, where Gallagher, Gallay, and Nier studied a one-dimensional simplified model 
\(-\partial_x^2 + x^2 + i\alpha f(x)\), approximating the operator under consideration by a complex harmonic oscillator.  
A similar strategy is frequently used in the semiclassical analysis of Schr\"odinger operators; see, for example, \cite{HS96}.

However, in the present setting we are dealing with the operator \(\mathcal{L}_1\), which has a complex-valued potential and is therefore non-self-adjoint. Consequently, the classical method of localizing the spectrum by constructing pseudo-eigenfunctions is no longer applicable.  
To overcome this difficulty, we introduce a new modification that remains effective in the non-self-adjoint framework, even in the presence of nonlocal terms (see Subsection~\ref{Spectral lower bound of k mode}).

We construct the following complex harmonic oscillator-type operator:
\begin{equation}\label{def:Z1}
\begin{aligned}
    \mathcal{Z}_1 u 
    \eqdefa &\Bigl(\partial_r^2 - \frac{35}{4r^2} - \frac{r^2}{16} + \frac{1}{2}\Bigr) u 
    - i \beta_1 \Bigl(\sigma(0) + \frac{\sigma''(0)}{2} r^2 \Bigr) u \\
    =& \Bigl(\partial_r^2 - \frac{35}{4r^2} - \frac{r^2}{16} + \frac{1}{2}\Bigr) u 
    - i \beta_1 \Bigl(1 - \frac{r^2}{8}\Bigr) u.
\end{aligned}
\end{equation}
In addition, the unitary scaling transformation is given by 
\begin{equation}
\label{unitary transformation}    (U_z u)(r) \eqdefa z^{1/2} u(zr), \quad z \in \R^+.
\end{equation}
In particular, the operators \(U_z\) are unitary on \(L^2_r\) and satisfy \(U_z^{-1} = U_{z^{-1}}\).  

We introduce the family of operators \(\{\mathcal{Z}_{1}^z\}_{z \in \R^+}\), which are mutually unitarily equivalent and therefore share the same spectrum, defined by  
\begin{equation}\label{def:Z1alpha_z}
    \mathcal{Z}_1^z u \eqdefa U_z \mathcal{Z}_1 U_z^{-1} u
    = z^{-2}\Bigl(\partial_r^2 - \frac{35}{4r^2}\Bigr)u
      - z^2\Bigl(\frac{1}{16} - \frac{i\beta_1}{8}\Bigr) r^2 u
      - i\beta_1 u + \frac12 u.
\end{equation}
The family \(\{\mathcal{Z}_{1}^z\}_{z \in \R^+}\) is then extended by analytic continuation to the sector  
\[
    \mathcal{S} \eqdefa \bigl\{ r e^{i\theta} : -\tfrac{\pi}{8} < \theta < \tfrac{\pi}{8},\; r > 0 \bigr\},
\]
and the resulting family is again denoted by \(\{\mathcal{Z}_{1}^z\}_{z \in \mathcal{S}}\).

One readily verifies that \(\{\mathcal{Z}_{1}^z\}_{z \in \mathcal{S}}\) forms an analytic family of type (A) in the sense of Kato \cite{Kato}, 
with common domain  
\[
    D_c \eqdefa \bigl\{ u \in H^2(\R^+; dr) :\; r^{-2}u,\; r^2 u \in L^2_r \bigr\}.
\]
Hence, by Kato’s theory, the spectrum of \(\mathcal{Z}_{1}^z\) is always discrete and depends holomorphically on \(z\).  
Since the eigenvalues of \(\mathcal{Z}_{1}^z\) remain constant for \(z \in \R^+\), they remain constant for all \(z \in \mathcal{S}\).

Moreover, since both \(f(r)\) and \(\sigma(r)\) admit natural analytic extensions from \(\R^+\) 
to the sector \(\mathcal{S}\), the same conclusions apply to the analytic continuation of 
the family \(\{\mathcal{L}_{1}^z\}_{z \in \R^+}\) to \(\{\mathcal{L}_{1}^z\}_{z \in \mathcal{S}}\), 
with the same domain \(D_c\). In particular,
\begin{equation}\label{def:L1_z}
    \mathcal{L}_{1}^z u
    \eqdefa U_z \mathcal{L}_1 U_z^{-1} u
    = z^{-2} \partial_r^2 u
      - \Bigl( z^{-2} \tfrac{3}{4r^2} + z^2 \tfrac{r^2}{16} - \tfrac{1}{2} + f(zr) \Bigr) u
      - i \beta_1 \sigma(zr) u.
\end{equation}
For additional details, the reader is referred to Section~7.1 of \cite{LWZ} 
and Appendix~\ref{basic calculus of complex deformation}.

We now fix 
\begin{equation}\label{def:zeta}
  \mathcal{S} \ni \zeta \eqdefa \Bigl(\frac{1}{16} - \frac{i\beta_1}{8}\Bigr)^{-\frac{1}{4}}
  \;\thicksim\; 2^{\frac{3}{4}} \beta_1^{-\frac{1}{4}} e^{\frac{\pi i}{8}} 
  + \cO(\beta_1^{-\frac{5}{4}}),
\end{equation}
and define
\begin{equation}\label{def:hatZ1}
     \hat{Z}_1 \eqdefa \zeta^2 \bigl(\mathcal{Z}_{1}^\zeta + i\beta_1 - \tfrac{1}{2}\bigr)
     = \partial_r^2 - \frac{35}{4r^2} - r^2,
\end{equation}
as well as
\begin{equation}\label{def:hatL}
   \hat{L}_1  \eqdefa \zeta^2 \bigl(\mathcal{L}_{1}^\zeta + i\beta_1 - \tfrac{1}{2}\bigr) 
     = \partial_r^2 - \Bigl( \tfrac{3}{4r^2} + \zeta^2 f(\zeta r) 
     + \zeta^4 \tfrac{r^2}{16} \Bigr) 
     + i\beta_1 \zeta^2 \bigl(1 - \sigma(\zeta r)\bigr).
\end{equation}

It is straightforward to verify that \(\hat{Z}_1\) is invertible, self-adjoint, negative, 
and possesses a compact resolvent. Its spectrum lies in \(\R^-\) and consists of a discrete 
unbounded set of eigenvalues with no accumulation point other than \(-\infty\). 
Consequently, the eigenvalues of \(\hat{Z}_1\) can be arranged in a decreasing sequence 
\(\{\lambda_k\}_{k \ge 1}\), and the multiplicity of \(\lambda_k\) is denoted by \(m_k\).

Our existence argument relies on the following fundamental lemma.
\begin{Lemma}\label{Lem 5.3}
    Let $\{A_\alpha\}_{\alpha > 0}$ be a family of compact operators on a Hilbert space $X$, 
    and let $A$ be a self-adjoint, negative, compact operator on $X$. 
    Let $\mu_k$ be the eigenvalues of $A$ in decreasing order,  with multiplicity $m_k$, and $d(\al)$ be any control function so that 
    \begin{equation*}
        \|A_\alpha - A\|_X \leq d(\al).
    \end{equation*}
    Then for any $0 < \delta \ll 1$ and all $\al$ satisfying  
    $d(\al) \leq \delta^2/4$, the following hold:
    \begin{itemize}
        \item[(i)] $\sigma(A_\alpha) \cap \Bigl( \bigcup_k B(\mu_k,\delta)\Bigr)^c = \emptyset$,
        \item[(ii)] for each open ball $B(\mu_k,\delta)$ satisfying  $\overline{B(\mu_k,\delta)} \cap \overline{B(\mu_j,\delta)} =\emptyset$ for $j\neq k$, the operator $A_\alpha$ 
        has eigenvalues inside $B(\mu_k,\delta)$ with total (algebraic) multiplicity $m_k$.
    \end{itemize}
\end{Lemma}
\begin{proof}
By Fredholm theory, both $\sigma(A_\alpha)$ and $\sigma(A)$ consist only of eigenvalues of finite multiplicity, 
with $0$ being the only possible accumulation point. 
Since $A$ is self-adjoint, we have for any $z \notin \sigma(A)$
\[
    \|(z - A)^{-1}\|_X^{-1} = \mathrm{dist}(z, \sigma(A)).
\]
If $d(\alpha) \leq \delta^2/4$ and 
$z \in \Bigl( \bigcup_k B(\mu_k, \delta)\Bigr)^c$, we obtain
\begin{equation}\label{eq:resolvent_estimate}
    \|(A - A_\alpha)(z - A)^{-1}\|_X 
    \leq \|A - A_\alpha\|_X \, \|(z - A)^{-1}\|_X
    \leq \frac{d(\alpha)}{\delta} 
    \leq \frac{\delta}{4}\ll1 ,
\end{equation}
which implies
\begin{equation*}
\begin{aligned}
\|(z - A_\alpha)^{-1}\|_X 
&= \bigl\| (z - A)^{-1} \bigl({\rm Id} +(A - A_\alpha)(z - A)^{-1}\bigr)^{-1}\bigr\|_X \\
&\leq  \|(z - A)^{-1}\|_X \times \bigl\|\bigl({\rm Id} + (A - A_\alpha)(z - A)^{-1}\bigr)^{-1}\bigr\|_X  \\
&\leq \frac{\|(z - A)^{-1}\|_X}{1 - \|(A - A_\alpha)(z - A)^{-1}\|_X} \leq \frac{1/\delta}{1 - \delta/4} \;\leq\; \frac{2}{\delta} \;<\; \infty.
\end{aligned}
\end{equation*}
This proves (i). 

To show (ii), we define  
\begin{equation}\label{def:riesz-projections}
    P \eqdefa \frac{1}{2\pi i} \int_{|z-\mu_k|= \delta} (z - A)^{-1} \, dz,
    \qquad 
    P_{\alpha} \eqdefa \frac{1}{2\pi i} \int_{|z-\mu_k|= \delta} (z - A_\alpha)^{-1} \, dz.
\end{equation}
It follows immediately from (i) that $|z-\mu_k|=\delta$ doesn't overlap with $\sigma(A_\al),\sigma(A)$,  for each open ball $B(\mu_k,\delta)$ satisfying  $\overline{B(\mu_k,\delta)} \cap \overline{B(\mu_j,\delta)} =\emptyset$ for $j\neq k$ and $d(\alpha) \leq \delta^2/4$. Thus
both $P$ and $P_\alpha$ are well-defined and represent the Riesz projections 
onto the corresponding root subspaces. In particular, we have 
$\dim \mathrm{Ran}(P) = m_k$ and
\begin{equation}\label{eq:projection-diff}
\begin{aligned}
    \|P_\alpha - P\|_X
    &\leq \max_{|z - \mu_k| = \delta} \|(z - A_\alpha)^{-1} - (z - A)^{-1}\|_X \\
    &\leq \max_{|z - \mu_k| = \delta} \|(z - A_\alpha)^{-1}\|_X \,\|A - A_\alpha\|_X \,\|(z - A)^{-1}\|_X \\
    &\leq \frac{2}{\delta} \cdot \frac{\delta^2}{4} \cdot \frac{1}{\delta} \;\leq\; \frac{1}{2}.
\end{aligned}
\end{equation}
In case  $\dim \mathrm{Ran}(P) \neq \dim \mathrm{Ran}(P_\alpha)$, 
without loss of generality, we may assume that $\dim \mathrm{Ran}(P) < \dim \mathrm{Ran}(P_\alpha) < \infty$. 
Then we can find some $\eta$ such that $P_\alpha \eta = \eta$ while $P \eta = 0$, 
which contradicts $\|P_\alpha - P\|_X \leq \tfrac{1}{2}$. 
This proves (ii) and we thus complete the proof of Lemma \ref{Lem 5.3}.
\end{proof}

The following coercive estimate indicates that the energy of the kernel function of $\hat{L}_1$, is mostly located in the region $r,s\thicksim\cO(1)$, will be crucial to our proof of Proposition \ref{Prop 5.3}.

\begin{Proposition}\label{Prop 5.4}
 There exists a constant $C > 0$, independent of $\alpha$, such that for $\alpha$ sufficiently large, the following bounds hold:
\begin{subequations}
    \begin{gather}
       \big\|\partial_r^2 \hat{Z}_1^{-1}\big\|_{\ltr}
    + \big\|r^2 \hat{Z}_1^{-1}\big\|_{\ltr}
    + \big\|r^{-2} \hat{Z}_1^{-1}\big\|_{\ltr}
    \;\leq\; C,   \label{eq 5.14a} \\ 
    \big\|\partial_r \hat{L}_1^{-1}\big\|_{\ltr}
    + \big\|r^{-1} \hat{L}_1^{-1}\big\|_{\ltr}
    + \big\|\hat{L}_1^{-1}\big\|_{\ltr}
    + \big\|\min\big\{r, \beta_1^{1/4}\big\} \, \hat{L}_1^{-1}\big\|_{\ltr}
    \;\leq\; C \label{eq 5.14b} \\
    \big\|\hat{L}_1^{-1} \min\big\{r, \beta_1^{1/4}\big\}\big\|_{\ltr}
    + \big\|\min\big\{r, \beta_1^{1/4}\big\} \, \hat{L}_1^{-1} \min\big\{r, \beta_1^{1/4}\big\} \big\|_{\ltr}
    \;\leq\; C. \label{eq 5.14c}
    \end{gather}
\end{subequations}
\end{Proposition}
\begin{proof}
In view of \eqref{def:hatZ1}, for any \(w \in D(\hat{Z}_1)\), taking the \(L^2_r\) inner product 
of \(\hat{Z}_1 w\) with \(r^2 w\) and using integration by parts yields
\begin{equation}\label{eq:Z1-coercive}
\begin{aligned}
\|\hat{Z}_1 w\|_{\ltr}\, \|r^2 w\|_{\ltr} 
&\;\geq\; -\Re\braket{\hat{Z}_1 w, r^2 w}_{\ltr} \\
&= \|r \partial_r w\|_{\ltr}^2 - \|w\|_{\ltr}^2 
   + 35/4\,\|w\|_{\ltr}^2 + \|r^2 w\|_{\ltr}^2 \\
&= \|r \partial_r w\|_{\ltr}^2 + 31/4\,\|w\|_{\ltr}^2 + \|r^2 w\|_{\ltr}^2,
\end{aligned}
\end{equation}
where \(\langle\cdot,\cdot\rangle_{\ltr}\) denotes the inner product in \(\ltr\).
In particular, this implies 
\[
    \|\hat{Z}_1 w\|_{\ltr} \ge \|r^2 w\|_{\ltr}, \quad \text{\rm and hence }\ 
   \big\| r^2 \hat{Z}_1^{-1} \big\|_{\ltr} \le 1.
\]
Proceeding in the same manner, taking the \(L^2_r\) inner product of \(\hat{Z}_1 w\) 
with \(r^{-2} w\) or with \(\partial_r^2 w\) yields the remaining parts of 
\eqref{eq 5.14a}. 

Next, in view of \eqref{def:hatL}, for any \(w \in D(\hat{L}_1)\),
\begin{equation*}
\begin{aligned}
\Re\braket{\hat{L}_1 w, w}_{\ltr}
&= - \|\partial_r w\|_{\ltr}^2 
   - 35/4\, \|r^{-1}w\|_{\ltr}^2
   - \braket{\Re\!\bigl(\zeta^2 f(\zeta r) - 8r^{-2}\bigr) w, w}_{\ltr} \\
&\quad - \Re(\zeta^4)/16\, \|r w\|_{\ltr}^2
   - \beta_1 \braket{\Im\!\bigl(\zeta^2(1 - \sigma(\zeta r))\bigr) w, w}_{\ltr}.
\end{aligned}
\end{equation*}
Since \(\Re(\zeta^4) > 0\), \(|\zeta| \sim \beta_1^{-1/4}\), and 
\(\zeta \in \mathcal{S}\) for \(\alpha \gg 1\), 
Lemmas~\ref{Lem A.1} and \ref{Lem A.2} give
\begin{equation*}
\begin{aligned}
 \|\partial_r w\|_{\ltr}^2
 + 35/4\,\|r^{-1}w\|_{\ltr}^2
 + \int_{\R^+} \min\{ r^2, \beta_1^{1/2} \}\, |w|^2\, dr
 \lesssim\; \bigl|\braket{\hat{L}_1 w, w}_{\ltr}\bigr| 
            + \beta_1^{-1/2}\|w\|_{\ltr}^2.
\end{aligned}
\end{equation*}
Combined with
\begin{equation*}
    \|r^{-1}w\|_{\ltr}^2
    + \int_{\R^+} \min\{ r^2, \beta_1^{1/2} \}\, |w|^2\, dr 
    \gtrsim \|w\|_{\ltr}^2,
\end{equation*}
it follows, for sufficiently large \(\alpha\), that
\begin{equation*}
\begin{split}
    &\|\partial_r w\|_{\ltr}^2 + \|w\|_{\ltr}^2
    + \|r^{-1}w\|_{\ltr}^2
    + \big\|\min\{ r, \beta_1^{1/4} \}\, w\big\|_{\ltr}^2 \\
    &\qquad\lesssim 
    \min\!\Big\{
        \|\hat{L}_1 w\|_{\ltr}\, \|w\|_{\ltr},\,
        \big\|\max\{ r^{-1}, \beta_1^{-1/4} \}\, \hat{L}_1 w\big\|_{\ltr}
        \big\|\min\{ r, \beta_1^{1/4} \}\, w\big\|_{\ltr}
    \Big\}.
\end{split}
\end{equation*}
This directly yields \eqref{eq 5.14b} and \eqref{eq 5.14c}, thereby completes the 
proof of Proposition~\ref{Prop 5.4}.
\end{proof}

\begin{Proposition}\label{Prop 5.5}
 There exists a constant \(C>0\), independent of \(\alpha\), such that for sufficiently large \(\alpha\),
\begin{equation}
    \big\| \hat{L}_1^{-1} - \hat{Z}_1^{-1} \big\|_{\ltr}
    \le C\, \beta_1^{-1/10}.
\end{equation}
\end{Proposition}
\begin{proof}
Smooth cut-off functions \(\chi_0\) and \(\chi_\infty\) are first chosen so that  
\(\chi_0 + \chi_\infty = 1\), with  
\(\chi_0 = \chi(r/\beta_1^{\kappa})\),  
where \(\chi\) is a smooth non-increasing function satisfying  
\(\chi = 1\) for \(r \in [0,1]\) and \(\chi = 0\) for \(r \in [2,\infty)\).  
The constant \(0 < \kappa < \tfrac14\) will be fixed later.  
We split
\begin{equation}\label{eq 5.17}
\begin{split}
\hat{L}_1^{-1} - \hat{Z}_1^{-1}
&= \hat{L}_1^{-1}(\hat{Z}_1 - \hat{L}_1)\hat{Z}_1^{-1} \\
&= \hat{L}_1^{-1}\chi_0(\hat{Z}_1 - \hat{L}_1)\hat{Z}_1^{-1}
 + \hat{L}_1^{-1}\chi_\infty(\hat{Z}_1 - \hat{L}_1)\hat{Z}_1^{-1}.
\end{split}
\end{equation}
Here,
\[
\hat{Z}_1 - \hat{L}_1
= \zeta^2 f(\zeta r) - \frac{8}{r^2}
  + i\beta_1 \zeta^2\!\left(\sigma(\zeta r) - 1 + \frac{\zeta^2 r^2}{8}\right)
\]
is purely a multiplication operator and will therefore be regarded as a function in the following.

Since \(|\zeta| \sim \beta_1^{-1/4}\), Proposition~\ref{Prop 5.4} together with  
Lemma~\ref{Lem A.1} yields
\begin{equation}\label{eq 5.18}
\begin{split}
\big\|\hat{L}_1^{-1}\chi_0(\hat{Z}_1-\hat{L}_1)\hat{Z}_1^{-1}\big\|_{\ltr}
&\lesssim \|\hat{L}_1^{-1}\|_{\ltr}\,
           \|\chi_0(\hat{Z}_1 - \hat{L}_1)\|_{\ltr}\,
           \|\hat{Z}_1^{-1}\|_{\ltr} \\
&\lesssim |\zeta|^2
 + \beta_1|\zeta|^2
    \min\!\big\{|\zeta|^4\beta_1^{4\kappa},\; |\zeta|^2\beta_1^{2\kappa}\big\} \\
&\lesssim \beta_1^{4\kappa - 1/2}.
\end{split}
\end{equation}

For the term involving \(\chi_\infty\), Proposition~\ref{Prop 5.4} gives
\begin{equation}\label{eq 5.19}
\begin{split}
\big\|\hat{L}_1^{-1}\chi_\infty(\hat{Z}_1-\hat{L}_1)\hat{Z}_1^{-1}\big\|_{\ltr}
&\lesssim \big\|\hat{L}_1^{-1}\chi_\infty\big\|_{\ltr}\,
           \Big\|\mathbbm{1}_{r\ge 1}(r/\beta_1^\kappa)
           (\hat{Z}_1 - \hat{L}_1)\hat{Z}_1^{-1}\Big\|_{\ltr} \\
&\lesssim \beta_1^{-\kappa}
           \big\|\hat{L}_1^{-1}\min\{r,\beta_1^{1/4}\}\chi_\infty\big\|_{\ltr}
           (|\zeta|^2 + \beta_1|\zeta|^4)\,
           \|r^2\hat{Z}_1^{-1}\|_{\ltr} \\
&\lesssim \beta_1^{-\kappa}.
\end{split}
\end{equation}

By combining \eqref{eq 5.18} and \eqref{eq 5.19} and choosing  
\(\kappa = \tfrac{1}{10}\), we achieve
\begin{equation*}
\|\hat{L}_1^{-1} - \hat{Z}_1^{-1}\|_{\ltr}
\lesssim \beta_1^{4\kappa - 1/2} + \beta_1^{-\kappa}
\lesssim \beta_1^{-1/10},
\end{equation*}
which completes the proof of Proposition~\ref{Prop 5.5}.
\end{proof}

We define the M{\"o}bius transformation $\Upsilon(z) = z^{-1}$ and then proceed to the proof of Proposition~\ref{Prop 5.3}.
\begin{proof}[\bf Proof of Proposition \ref{Prop 5.3}]
Let 
$$ \beta_1=\frac{\alpha}{8\pi}, \quad \zeta=\big(\frac{1}{16}-\frac{i\beta_1}{8}\big)^{-\frac{1}{4}}, \quad A_\alpha=\hat{L}_1^{-1},\quad A=\hat{Z}_1^{-1},\quad  d(\alpha)=C|\alpha|^{-\frac{1}{10}},$$
where thanks to Proposition \ref{Prop 5.5}, we can take  the constant $C$ to be independent of $\al$ and large enough so that $\|A_\al-A\|_{\ltr}\leq d(\al).$  We denote $0>\la_1>\la_2>\cdots\to -\infty$ to be the eigenvalues of $\hat{Z}_1=\partial_r^2 - \frac{35}{4r^2}-r^2$ with  multiplicities $m_k$. Thus for any $N>0$, there exists $\delta_N>0$  so that  $\overline{B(\la_k^{-1},\delta_N)} \cap \overline{B(\la_j^{-1},\delta_N)}=\emptyset$ for $1\leq k\leq N$ and all $j\neq k$. Without loss of generality, we always take this $\delta_N$ to be sufficiently small.

In order to obtain spectral properties of $\hat{L}_1$, and hence of $\hat{L}_1^{-1}$, we deduce from Lemma \ref{Lem 5.3} and 
\[
\sigma_{X_1\cap Y_1}(L_1-\alpha\Lambda_1)
= \sigma_{\ltr}(\mathcal{L}_{1})
= \sigma_{\ltr}(\mathcal{L}_{1}^\zeta)
= \zeta^{-2}\sigma_{\ltr}(\hat{L}_1)-i\beta_1+\tfrac{1}{2},
\]
that  for any $\alpha \geq (4C)^{10}\delta_N^{-20}$, the following hold:

\noindent{\bf (1).} $\sigma_{\ltr}(\mathcal{L}_{1}) \subset 
    \zeta^{-2}\!\bigcup_{k=1}^\infty \Upsilon\!\big(B(\la_k^{-1},\delta_N)\big)
    - i\beta_1 + \tfrac{1}{2};$
    
\noindent{\bf (2).} for $1\leq k\leq N$, there exist eigenvalues of $\cL_1$ of total (algebraic) multiplicity $m_k$ inside the set
$\zeta^{-2}\Upsilon\!\big(B(\la_k^{-1},\delta_N)\big) 
    - i\beta_1 + \frac{1}{2}.$
    
In particular, there exists
\begin{equation*}
    z^\star =\zeta^{-2}w^\star-i \beta_1+1/2 \with z^\star\in \sigma_{\ltr}(\cL_1) \andf  w^\star \in \Upsilon\!\big(B(\la_1^{-1},\delta_N)\big).
\end{equation*}
Notice that $\la_1 < 0 $, $\delta_N\ll1$ and recall \eqref{def:zeta},
we achieve 
\begin{equation*}
   s_{\ltr}(\cL_1)\geq {\rm Re} \, z^\star = {\rm Re}\, (\zeta^{-2}w^\star)+1/2 =  \beta_1^{1/2}/4 \times  {\rm Re}\, \big( (\beta_1^{-1}-2i)^{1/2} w^\star  \big) +1/2 \geq -C \beta_1^{1/2}.
\end{equation*}
This completes the proof of Proposition \ref{Prop 5.3}.
\end{proof}

Now we are in a position to present the proof of Theorem \ref{Thm 2}.
\begin{proof}[\bf Proof of Theorem \ref{Thm 2}]
Noting that  
\begin{equation*}
\begin{split}
\sigma_{({\rm Ker}\,\La)^\perp}(L_\perp - \alpha \La_\perp)
= \sigma_{X_1 \cap Y_1}(L_1 - \alpha \La_1)
  \cup \sigma_{X_{-1} \cap Y_1}(L_{-1} - \alpha \La_{-1})
  \cup_{k \neq \pm 1}\, \sigma_{X_k}(L_k - \alpha \La_k),
\end{split}
\end{equation*}
Proposition~\ref{Prop 5.3} implies that  
\begin{equation*}
\Sigma(\alpha)
= - s_{({\rm Ker}\,\La)^\perp}(L_\perp - \alpha \La_\perp)
\le - s_{X_1 \cap Y_1}(L_1 - \alpha \La_1)
\le C |\alpha|^{1/2}.
\end{equation*}
Combined with Theorem~\ref{Thm 1.3}, this concludes the proof of Theorem~\ref{Thm 2}.
\end{proof}

\begin{figure}[H]
 \centering
 \includegraphics[width=9cm]{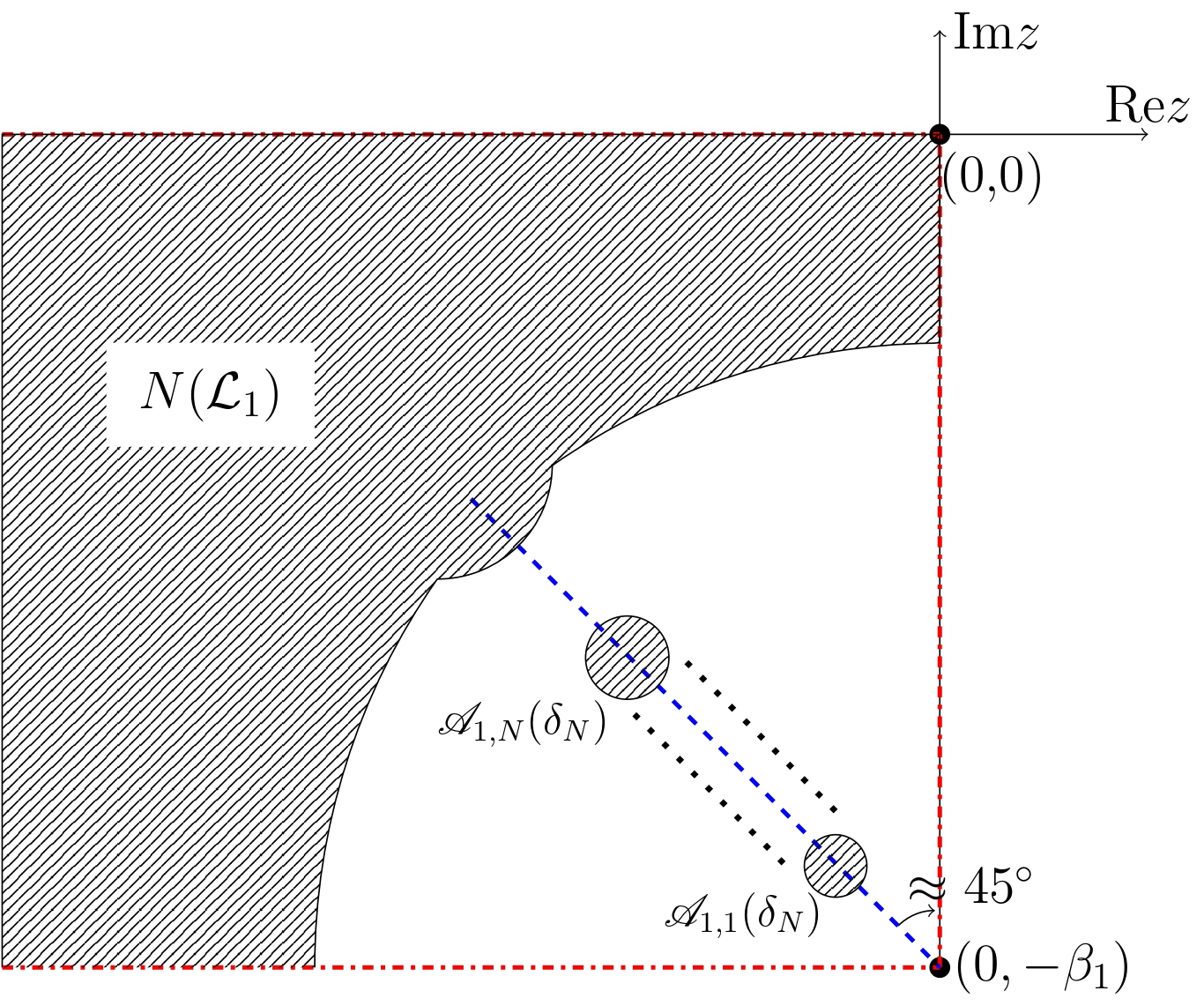}
 \caption{\small Spectral localization of $\sigma_{\ltr}(\cL_1)$: the spectrum $\sigma_{\ltr}(\cL_1)$ lies entirely within the shaded region, and each connected component of this region contains at least one spectral point of $\sigma_{\ltr}(\cL_1)$.} 
\label{fig:spectrum1}
 \end{figure}

We end this subsection with several interesting remarks:
\begin{itemize}
    \item[(1)]Observe that 
\[
U_2^{-1} \hat{Z}_1 U_2
= 4\Bigl(\partial_r^2 - \frac{35}{4r^2} - \frac{r^2}{16}\Bigr)
= 4 r^{1/2} e^{r^2/8}\, L_3 \, r^{-1/2} e^{-r^2/8} - 2.
\]
If \(f(r)\) is an eigenfunction of \(L_k\), then  
\(f(r)e^{ik\theta}\) must be an eigenfunction of \(L\), and therefore must be of the form  
\(\nabla_\xi^\beta G\) for some \(\beta \in \mathbb{N}^2\).
This implies that \(f(r)\) necessarily has the structure  
\(r^{|k|} Q_{\ell,k}(r)e^{-r^2/4}\).
Substituting this into the equation \(L_k f = \lambda f\) allows one to compute \((\lambda,f(r))\) explicitly.

For example, the functions  
\(r^{|k|} e^{-r^2/4}\) and  
\(r^{|k|}(r^2 - 4|k| - 4)e^{-r^2/4}\)  
are, by Sturm--Liouville theory, eigenfunctions of \(L_k\) corresponding to the first two eigenvalues  
\(-\frac{|k|}{2}\) and \(-\frac{|k|+2}{2}\), respectively.
Hence, using the spectral properties of \(L_k\), one obtains \(\lambda_1 = -8\). The remaining eigenvalues \(\lambda_k\) of \(\hat{Z}_1\), together with their corresponding eigenfunctions, can also be computed explicitly, although these computations will not be pursued here.
\item[(2)]
Moreover, it is known that  
\begin{equation}\label{def:numerical range}
\begin{split}
     N(\mathcal{L}_{1})
     &\eqdefa \bigl\{ \braket{\mathcal{L}_{1}u, u}_{\ltr} 
         :\ u \in D(\mathcal{L}_{1}),\ \|u\|_{\ltr}=1 \bigr\} \\
     &\subset \{ z \in \mathbb{C} :\ \Re z \le 0,\ -\beta_1 \le \Im z \le 0 \},
\end{split}
\end{equation}
which implies  
\begin{equation*}
    \sigma_{\ltr}(\mathcal{L}_{1})
        \subset \{ z \in \mathbb{C} :\ \Re z \le 0,\ -\beta_1 \le \Im z \le 0 \}.
\end{equation*}
Let  
\[
\mathscr{A}_{1,j}(\delta_N)
\eqdefa \zeta^{-2}\Upsilon\!\bigl(B(\lambda_j^{-1},\delta_N)\bigr)
        - i\beta_1 + \tfrac12,
        \qquad j \le N.
\]
A more refined localization of \(\sigma_{\ltr}(\mathcal{L}_1)\) is then obtained;
see Figure~\ref{fig:spectrum1}.
\end{itemize}

\subsection{Spectral lower bound of $L_k - \al \Lambda_k\ (k\geq2)$}\label{Spectral lower bound of k mode} 
For the proof of Theorem~\ref{Thm 2}, the results obtained in Subsection~\ref{sect. 5.1} 
are already sufficient. Nevertheless, in this subsection a similar result is established 
for \(H_k\). It should be emphasized that, in this case, \(H_k\) contains a genuinely 
nonlocal term. Thus, the argument provided here serves both for completeness and for its 
own independent interest, while also illustrating the robustness of the method.

\begin{Proposition}\label{Prop 5.6}
 There exists a constant $C > 0$ such that for any $k \geq 2$,
\[
s_{X_k}(L_k - \alpha \Lambda_k) = s_{L^2_r}(H_k) \geq - C |k\alpha|^{1/2}.
\]
\end{Proposition}
Recall that $\beta_k = \frac{k\alpha}{8\pi}$ and
\begin{equation*}
    H_{k} 
    = \partial_r^2 - \Big(k^2-\frac{1}{4}\Big)r^{-2} - \frac{r^2}{16} + \frac{1}{2} 
    - i\beta_k \bigl(\sigma(r) - g(r)\mathcal{K}_k\big[g(r)\cdot\big]\bigr).
\end{equation*}
Following the same approach as in the previous subsection, we construct
\begin{equation*}
    \mathcal{Z}_{k} 
    \eqdefa \partial_r^2 -  \Big(k^2-\frac{1}{4}\Big)r^{-2} - \frac{r^2}{16} + \frac{1}{2}
    - i\beta_k \Bigl( 1 - \frac{r^2}{8} - \mathcal{K}_k \Bigr).
\end{equation*} 

Similarly, we have
\begin{equation*}
\begin{split}
    \mathcal{Z}_{k}^z 
    &\eqdefa U_z \, \mathcal{Z}_{k} \, U_{z}^{-1} \\
    &= z^{-2}\Bigl(\partial_r^2 -  \Big(k^2-\frac{1}{4}\Big)r^{-2}\Bigr) 
       - z^2\Bigl(\frac{1}{16} - \frac{i\beta_k}{8}\Bigr)r^2 
       + i\beta_k z^2 \mathcal{K}_k - i\beta_k + \frac{1}{2},
\end{split}
\end{equation*}
and
\begin{equation*}
\begin{split}
    H_{k}^z 
    &\eqdefa U_z \, H_{k} \, U_{z}^{-1} \\
    &= z^{-2}\Bigl(\partial_r^2 - \Big(k^2-\frac{1}{4}\Big)r^{-2}\Bigr) 
       - \frac{z^2 r^2}{16} + \frac{1}{2} 
       - i\beta_k \bigl(\sigma(zr) - z^2 g(zr)\mathcal{K}_k\big[g(zr)\cdot\big]\bigr),
\end{split}
\end{equation*}
which form analytic families of type (A) for $z \in \mathcal{S}$, with common domain $D_c$ (cf.~Kato~\cite{Kato}). 
Therefore, the spectra of $\{\mathcal{Z}_{k}^z\}_{z \in \mathcal{S}}$ and $\{H_{k}^z\}_{z \in \mathcal{S}}$ consist only of discrete eigenvalues and are independent of $z \in \mathcal{S}$. Now we fix  
\[
\cS \ni \zeta_k \eqdefa \Bigl(\frac{1}{16} - \frac{i\beta_k}{8}\Bigr)^{-\frac{1}{4}} \thicksim 2^{\f34}\beta_k^{-\f14}e^{\frac{i\pi }{8}} + \cO(\beta_k^{-\f54}),
\]  
so that $\mathcal{Z}_{k}^{\zeta_k}$ takes the form  
\begin{equation*}
\mathcal{Z}_{k}^{\zeta_k} 
= \zeta_k^{-2}\Bigl(\partial_r^2 - \Big(k^2-\frac{1}{4}\Big)r^{-2} - r^2 
  + \frac{-32\beta_k^2+16i\beta_k}{1+4\beta_k^2}\mathcal{K}_k\Bigr) 
  - i\beta_k + \frac{1}{2}.
\end{equation*}

Finally, we construct the operator  
\begin{equation}
\label{Zk}\hat{Z}_k \eqdefa \partial_r^2 - \Big(k^2-\frac{1}{4}\Big)r^{-2} - r^2 - 8\mathcal{K}_k,
\end{equation}
which serves as an approximate model for the operator $\zeta_k^{2}(\mathcal{Z}_{k}^{\zeta_k}+i\beta_k-\frac{1}{2})$. In parallel, we introduce  
\begin{equation*}
\begin{split}
\hat{L}_k 
&\eqdefa \zeta_k^2 \bigl(H_{k}^{\zeta_k} + i\beta_k - 1/2\bigr) \\
&= \partial_r^2 - \Big(k^2-\frac{1}{4}\Big)r^{-2} -  \frac{\zeta_k^4}{16}r^2 
   + i\beta_k \zeta_k^2\bigl(1-\sigma(\zeta_k r)\bigr) 
   + i\beta_k \zeta_k^4 g(\zeta_k r)\mathcal{K}_k[g(\zeta_k r)\cdot].
\end{split}
\end{equation*}

The following basic proposition establishes the boundedness and coercivity of the inverses $\hat{L}_k^{-1}$ and $\hat{Z}_k^{-1}$.
\begin{Proposition}\label{Prop 5.7}
 For all $|k| \ge 2$, $\exists$ a constant $C>0$, independent of $k$ and $\alpha$, such that
\begin{subequations}
    \begin{gather}
        k^2 \|\mathcal{K}_k f\|_{\ltr} \le C \|r^2 f\|_{\ltr}, \label{eq 5.21a} \\
        \big\|\partial_r^2 \hat{Z}_k^{-1}\big\|_{\ltr}
        + k^2 \bigl\|r^{-2}\hat{Z}_k^{-1}\bigr\|_{\ltr}
        + \big\|r^2 \hat{Z}_k^{-1}\big\|_{\ltr}
        + |k|\,\big\|\hat{Z}_k^{-1}\big\|_{\ltr}
        \le C, \label{eq 5.21b} \\
        \big\|\partial_r \hat{L}_k^{-1}\big\|_{\ltr}
        + \bigl\|r^{-1}\hat{L}_k^{-1}\bigr\|_{\ltr}
        + \big\|\hat{L}_k^{-1}\big\|_{\ltr}
        + \big\|\min\{r,\beta_k^{1/4}\}\,\hat{L}_k^{-1}\big\|_{\ltr}
        \le C, \label{eq 5.21c} \\
        \big\|\hat{L}_1^{-1}\,\min\{r,\beta_k^{1/4}\}\big\|_{\ltr}
        + \big\|\min\{r,\beta_k^{1/4}\}\,\hat{L}_1^{-1}\,\min\{r,\beta_k^{1/4}\}\big\|_{\ltr}
        \le C. \label{eq 5.21d}
    \end{gather}
\end{subequations}
\end{Proposition}
\begin{proof} 
The identity $\widetilde{\Delta}_k = \partial_r^2 - (k^2-\tfrac14) r^{-2}$ is first recalled.  
For $|k|\ge 2$ and $f\in D(\widetilde{\Delta}_k)$, 
\begin{equation}\label{estimate of lap_k}
\begin{split}
  \|r^2 \widetilde{\Delta}_k f\|_{\ltr} \, \|f\|_{\ltr} 
  &\ge -{\rm Re}\,\braket{\widetilde{\Delta}_k f, r^2 f}_{\ltr} 
   = \|r \partial_r f\|_{\ltr}^2 + (k^2-5/4)\|f\|_{\ltr}^2 .
\end{split}
\end{equation}
Since $-\widetilde{\Delta}_k\mathcal{K}_k = {\rm Id}$, the estimate \eqref{eq 5.21a} follows directly from \eqref{estimate of lap_k}.

To obtain \eqref{eq 5.21b}, observe that for $f\in D_c$,
\begin{equation*}
\begin{split}
\|\hat{Z}_k f\|_{\ltr}\, \|\widetilde{\Delta}_k f\|_{\ltr}
&\ge {\rm Re}\,\braket{\hat{Z}_k f, \widetilde{\Delta}_k f}_{\ltr} \\
&= \|\widetilde{\Delta}_k f\|_{\ltr}^2
   + \|r\partial_r f\|_{\ltr}^2
   + (k^2-5/4)\|f\|_{\ltr}^2
   + 8\|f\|_{\ltr}^2 .
\end{split}
\end{equation*}
Thanks to 
\[
\|\widetilde{\Delta}_k f\|_{\ltr}
\thicksim \|\partial_r^2 f\|_{\ltr}
+ |k|\,\|r^{-1}\partial_r f\|_{\ltr}
+ k^2\|r^{-2} f\|_{\ltr},
\qquad |k|\ge 2,
\]
one obtains \eqref{eq 5.21b} except for the third term, which follows in the same manner by taking the $L^2_r$ inner product of $\hat{Z}_k f$ with $r^2 f$ and repeating the above argument.

The estimate \eqref{eq 5.21c} is derived next.  
For any $f\in D_c$,
\begin{equation*}
\begin{split}
-{\rm Re}\,\braket{\hat{L}_k f, f}_{\ltr}
&= \|\partial_r f\|_{\ltr}^2
   + (k^2-1/4)\|r^{-1}f\|_{\ltr}^2
   + \frac{{\rm Re}\,\zeta_k^4}{16}\|r f\|_{\ltr}^2 \\
&\quad + \beta_k \int_0^\infty 
      {\rm Im}\bigl(\zeta_k^2(1-\sigma(\zeta_k r))\bigr)\, |f(r)|^2\, dr \\
&\quad + \beta_k \int_0^\infty\!\!\int_0^\infty 
      K_k(r,s)\,{\rm Im}\bigl(\zeta_k^4 g(\zeta_k r) g(\zeta_k s)\bigr)\,
      f(r)\,\overline{f(s)}\, dr ds .
\end{split}
\end{equation*}
For sufficiently large $\alpha$, $\zeta_k$ lies in $\mathcal{S}_4$.  
Inserting \eqref{ineq, coercivity with nonlocal term} into the above inequality and using ${\rm Re}\,\zeta_k^4>0$ gives
\begin{equation*}
\begin{split}
\big|{\rm Re}\,\braket{\hat{L}_k f,f}_{\ltr}\big|
&\gtrsim \|\partial_r f\|_{\ltr}^2
   + (k^2-1/4)\|r^{-1}f\|_{\ltr}^2
   + \beta_k \int_0^\infty \min\{|\zeta_k|^4 r^2, |\zeta_k|^2\}\, |f(r)|^2\, dr \\
&\gtrsim \|\partial_r f\|_{\ltr}^2
   + (k^2-1/4)\|r^{-1}f\|_{\ltr}^2
   + \int_0^\infty \min\{r^2, \beta_k^{1/2}\}\, |f(r)|^2\, dr .
\end{split}
\end{equation*}
Estimates \eqref{eq 5.21c} and \eqref{eq 5.21d} then follow by the same lines as in the proof of Proposition~\ref{Prop 5.5}.  
This completes the proof of Proposition~\ref{Prop 5.7}.
\end{proof}

\begin{Proposition}\label{Prop 5.8}
For $\alpha$ sufficiently large, there exists a constant $C>0$, independent of $k$ and $\alpha$, such that
\[
   \big\|\hat{L}_k^{-1}-\hat{Z}_k^{-1}\big\|_{\ltr} \le C\, \beta_k^{-1/10}.
\]
\end{Proposition}
\begin{proof}
Let $0<\kappa < \frac{1}{4}$. Define $\chi_0$ and $\chi_\infty$ in the same manner as in the proof of Proposition~\ref{Prop 5.5}, with $\beta_1$ replaced by $\beta_k$. Observing that 
\[
     \hat{L}_k^{-1}-\hat{Z}_k^{-1} = \hat{L}_k^{-1}(\hat{Z}_k-\hat{L}_k)\hat{Z}_k^{-1},
\]
and
\begin{equation*}
\begin{split}
    \hat{Z}_k-\hat{L}_k
    &= i\beta_k\zeta_k^2\bigl(\sigma(\zeta_k r)-1+\zeta_k^2 r^2/8\bigr) 
       -(8+i\beta_k\zeta_k^4)g(\zeta_k r)\mathcal{K}_k[g(\zeta_k r)\cdot] \\
    &\quad + 8\bigl(g(\zeta_k r)-1\bigr)\mathcal{K}_k[g(\zeta_k r)\cdot] 
       + 8\mathcal{K}_k[(g(\zeta_k r)-1)\cdot] \eqdefa A_1+\cdots+A_4,
\end{split}
\end{equation*}
we split
\begin{equation*}
    \begin{split}
        \hat{L}_k^{-1}-\hat{Z}_k^{-1}
        &= \hat{L}_k^{-1}\chi_0 A_1\hat{Z}_k^{-1}
         + \cdots
         + \hat{L}_k^{-1}\chi_0 A_4\hat{Z}_k^{-1} 
         + \hat{L}_k^{-1}\chi_\infty (\hat{Z}_k-\hat{L}_k)\hat{Z}_k^{-1} \\
        &\eqdefa O_1+O_2+O_3+O_4+O_5.
    \end{split}
\end{equation*}

Using Proposition~\ref{Prop 5.7} together with Lemma~\ref{Lem A.1}, one obtains
\begin{equation}\label{eq 5.23}
\begin{split}
    \|O_1\|_{\ltr}
    &\lesssim \beta_k |\zeta_k|^2 
       \min\big\{|\zeta_k|^4 \beta_k^{4\kappa},\; |\zeta_k|^2 \beta_k^{2\kappa}\big\}
       \|\hat{L}_k^{-1}\|_{\ltr}\,\|\hat{Z}_k^{-1}\|_{\ltr}
       \lesssim \beta_k^{4\kappa-1/2}.
\end{split}
\end{equation}

For $O_2$, using $|8+i\beta_k\zeta_k^4|\lesssim \beta_k^{-1}$, $|g(\zeta_k r)|\lesssim 1$, and Proposition~\ref{Prop 5.7}, we deduce
\begin{equation}
\begin{split}
    \|O_2\|_{\ltr}
    &\lesssim \beta_k^{-1}\|\hat{L}_k^{-1}\|_{\ltr}
       \big\|g(\zeta_k r)\mathcal{K}_k[g(\zeta_k r)\hat{Z}_k^{-1}\cdot]\big\|_{\ltr} \\
    &\lesssim \beta_k^{-1}\|\hat{L}_k^{-1}\|_{\ltr}\,\|r^2\hat{Z}_k^{-1}\|_{\ltr}
       \lesssim \beta_k^{-1}.
\end{split}
\end{equation}

For $O_3$, using $|g(\zeta_k r)-1|\lesssim \min\{|\zeta_k|^2 r^2,1\}$, we similarly obtain
\begin{equation}
\begin{split}
    \|O_3\|_{\ltr}
    &\lesssim |\zeta_k|^2\beta_k^{2\kappa}
       \|\hat{L}_k^{-1}\|_{\ltr}\,
       \|\mathcal{K}_k[g(\zeta_k r)\hat{Z}_k^{-1}\cdot]\|_{\ltr} \lesssim \beta_k^{2\kappa-1/2}.
\end{split}
\end{equation}

For $O_5$, Proposition~\ref{Prop 5.7} and Lemma~\ref{Lem A.1} imply
\begin{equation}
\begin{split}
    \|O_5\|_{\ltr}
    &\lesssim \|\hat{L}_k^{-1}\chi_\infty\|_{\ltr}
       \big\|\mathbbm{1}_{r\ge 1}(r/\beta_k^{\kappa})(\hat{Z}_k-\hat{L}_k)\hat{Z}_k^{-1}\big\|_{\ltr} \\
    &\lesssim \beta_k^{-\kappa}
       \|\hat{L}_1^{-1}\min\{r, \beta_k^{1/4}\}\|_{\ltr}\,
       \|r^2\hat{Z}_k^{-1}\|_{\ltr}
       \lesssim \beta_k^{-\kappa}.
\end{split}
\end{equation}

To estimate $O_4$, we first introduce $\iota_0$ and $\iota_\infty$ in the same way as $\chi_0$ and $\chi_\infty$, but with $\kappa$ replaced by $2\kappa$. Then
\[
O_4
= \hat{L}_k^{-1}\chi_0\mathcal{K}_k[\iota_0(g(\zeta_k r)-1)\hat{Z}_k^{-1}\cdot]
 + \hat{L}_k^{-1}\chi_0\mathcal{K}_k[\iota_\infty(g(\zeta_k r)-1)\hat{Z}_k^{-1}\cdot]
 \eqdefa O_{41}+O_{42}.
\]
As in the estimate of $O_3$, we find
\begin{equation}
\begin{split}
    \|O_{41}\|_{\ltr}
    &\lesssim \|\hat{L}_k^{-1}\|_{\ltr}\,
             \|\mathcal{K}_k[\iota_0(g(\zeta_k r)-1)\hat{Z}_k^{-1}\cdot]\|_{\ltr} \\
    &\lesssim \|\hat{L}_k^{-1}\|_{\ltr}\,
             \|r^2\iota_0(g(\zeta_k r)-1)\hat{Z}_k^{-1}\|_{\ltr} \\
    &\lesssim \|\hat{L}_k^{-1}\|_{\ltr}\,
             \|\iota_0(g(\zeta_k r)-1)\|_{\ltr}\,
             \|r^2\hat{Z}_k^{-1}\|_{\ltr}
             \lesssim \beta_k^{4\kappa-1/2}.
\end{split}
\end{equation}

Using the explicit kernel representation of $\mathcal{K}_k$, we obtain
\begin{equation*}
\begin{split}
\big|\chi_0\mathcal{K}_k[\iota_\infty(g(\zeta_k r)-1)\hat{Z}_k^{-1}w]\big|(r)
&\lesssim \chi_0(r)\, r^{5/2}
   \int_{\beta_k^{2\kappa}}^\infty s^{-7/2} s^2 |\hat{Z}_k^{-1}w|(s)\, ds \\
&\lesssim \chi_0(r)\, r^{5/2}\, \beta_k^{-6\kappa}\, \|r^2\hat{Z}_k^{-1}w\|_{\ltr}.
\end{split}
\end{equation*}
Consequently,
\begin{equation}\label{eq 5.28}
\begin{split}
    \|O_{42}w\|_{\ltr}
    &\lesssim \beta_k^{-6\kappa}\,
               \|\hat{L}_k^{-1}\|_{\ltr}\,
               \|\chi_0 r^{5/2}\|_{\ltr}\,
               \|r^2\hat{Z}_k^{-1}w\|_{\ltr} \lesssim \beta_k^{-3\kappa}\,\|w\|_{\ltr}.
\end{split}
\end{equation}

Combining \eqref{eq 5.23}--\eqref{eq 5.28}, we conclude that
\[
    \|\hat{L}_k^{-1}-\hat{Z}_k^{-1}\|_{\ltr}
    \lesssim \beta_k^{4\kappa-1/2} + \beta_k^{-\kappa},
\]
and choosing $\kappa=\frac{1}{10}$ completes the proof of Proposition~\ref{Prop 5.8}.
\end{proof}

\medskip

\begin{figure}[H]
\centering
\includegraphics[width=9cm]{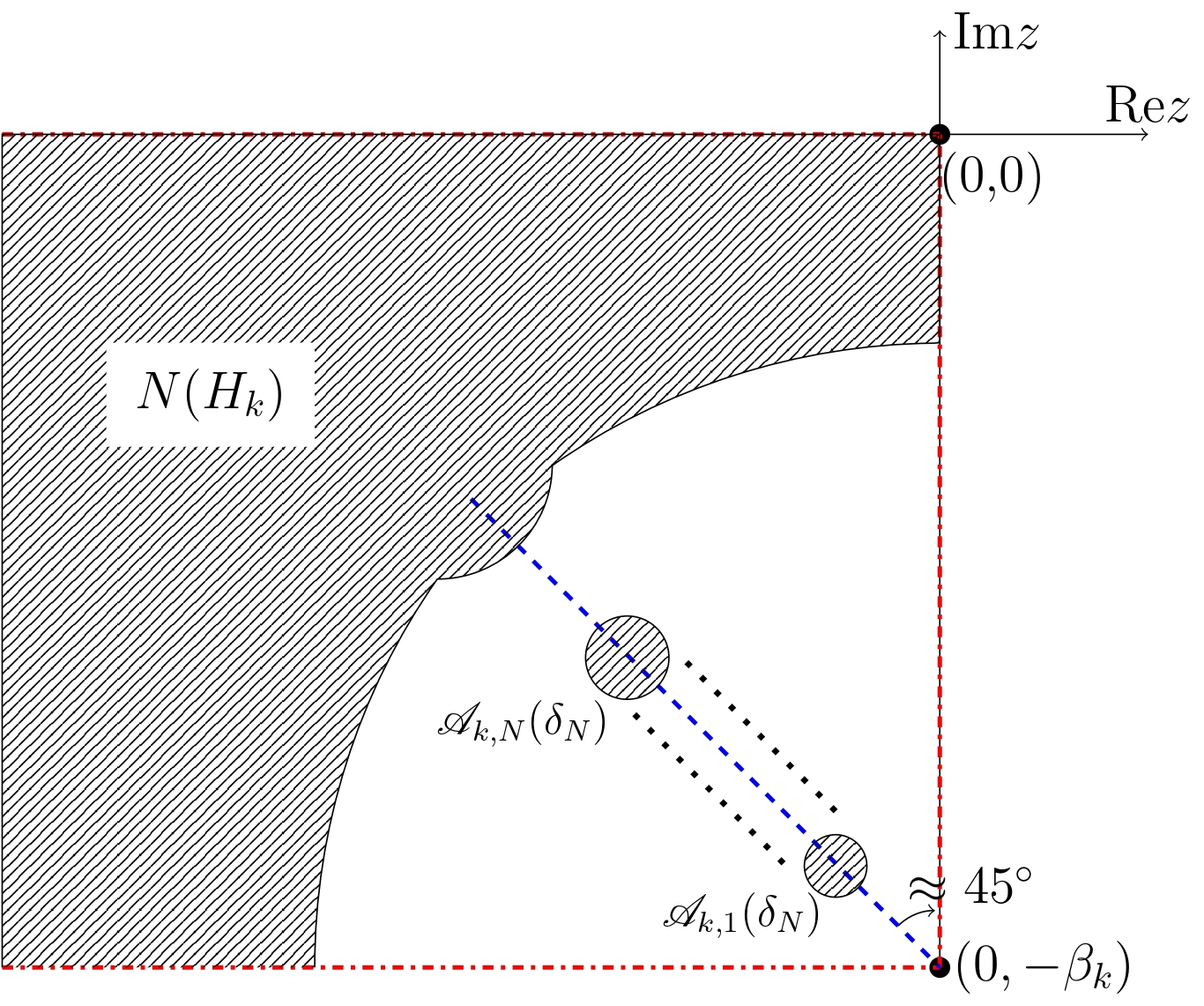}
\caption{Spectral localization of $\sigma_{\ltr}(H_{k})$: the spectrum $\sigma_{\ltr}(H_{k})$ lies within the shaded region, and each connected component of this region contains at least one spectral point.} 
\label{fig:Hk-spectrum}
\end{figure}

Let us now present the proof of Proposition \ref{Prop 5.6}.
\begin{proof}[\bf Proof of Proposition \ref{Prop 5.6}]
Let 
\[
\beta_k=\frac{k\alpha}{8\pi}, \quad 
\zeta_k=\left(\frac{1}{16}-\frac{i\beta_k}{8}\right)^{-\frac{1}{4}}, \quad 
A_\alpha=\hat{L}_k^{-1},\quad 
A=\hat{Z}_k^{-1},\quad  
d(\alpha)=C|\beta_k|^{-\frac{1}{10}},
\]
where the constant $C$ is independent of $\alpha$ and $k$, and is chosen sufficiently large so that $\|A_\alpha-A\|_{\ltr}\le d(\alpha)$ by Proposition~\ref{Prop 5.8}.  
Let $0>\tilde{\lambda}_1>\tilde{\lambda}_2>\cdots\to -\infty$ denote the eigenvalues of
\[
\hat{Z}_k=\partial_r^2 - \frac{4k^2-1}{4r^2} - r^2 - 8\mathcal{K}_k,
\]
with multiplicities $m_k$.  
The proof of Proposition~\ref{Prop 5.6} then follows the same lines as the proof of Proposition~\ref{Prop 5.3}.  
For $j\le N$, define
\[
\mathscr{A}_{k,j}(\delta_N)\eqdefa 
\zeta_k^{-2}\Upsilon\!\big(B(\tilde{\lambda}_j^{-1},\delta_N)\big)
    - i\beta_k + \frac{1}{2},
\]
and let $N(H_k)$ denote the numerical range of $H_k$, defined in the same way as in \eqref{def:numerical range}.  
The location of $\sigma_{\ltr}(H_k)$ is shown in Figure~\ref{fig:Hk-spectrum}.
\end{proof}

\section{Spectral bound of $\sL-\al\Pi$.}
In this section, we investigate the spectral bound of $\sL-\al\Pi$. In particular, we establish the following estimates.
\begin{Lemma}\label{Lem 3.2}
 There exists a constant $C>0$, independent of $k$ and $\al$, such that
\begin{equation}\label{eq 3.2}
    s_{X_k^2}(\sL_k-\al \Pi_k) \leq -C |k\al|^{1/2}.
\end{equation}
Consequently,
\begin{equation}\label{eq 3.3}
    s_{X_{\neq}^2}(\sL-\al \Pi) \leq -C |\al|^{1/2}.
\end{equation}
\end{Lemma}

\begin{proof}
We first deduce, by compact perturbation theory together with the spectral properties of $L$ and the definition of $\sL$, that $\sigma_{Y^2}(\sL-\al \Pi)$ consists solely of eigenvalues with finite multiplicity.  
Thus it suffices to prove \eqref{eq 3.2}.  

Let $\la \in \sigma_{X_k^2}(\sL_k-\al\Pi_k)$.  
Then there exists a nonzero eigenfunction $f=(f_1,f_2)^\tr \in D(\sL_k-\al\Pi_k)\cap X_{\neq}^2$ satisfying  
\[
(\sL_k-\al\Pi_k)f=\la f.
\]
The observation in \cite{GM11, BGH23} shows that
\[
f_{\rm div} \eqdefa r^{-1}\pa_r(rf_1) + i k r^{-1} f_2
\]
satisfies
\begin{equation}\label{eq B.3}
    (L_k + 1/2 - i \al k S(r)) f_{\rm div} = \la f_{\rm div}.
\end{equation}
Together with Lemma~\ref{Lem A.4}, this yields \eqref{eq 3.2}, which completes the proof of Lemma~\ref{Lem 3.2}.
\end{proof}

Following the strategy in \cite{LWZ}, we establish the next lemma.
\begin{Lemma}\label{Lem A.4}
 Let $|\al| \gg 1$. Then there exists a universal constant $C>0$, independent of $\al$ and $k$, such that
\begin{equation}
    s_Y(L_k - i\al k S(r)) \leq -C |\al k|^{1/2}.
\end{equation}
\end{Lemma}

\begin{proof}
It suffices to consider the action of $L_k - i\al k S(r)$ on the space $L^2(\R^+; r e^{r^2/4}\,dr)$, or equivalently, the operator
\begin{equation*}
    \mathscr{H}_{\al,k} \eqdefa 
    e^{r^2/8} r^{1/2} (L_k - i\al k S(r)) r^{-1/2} e^{-r^2/8}
    = \pa_r^2 - \Big(k^2-\tfrac14\Big) r^{-2} - \frac{r^2}{16} + \frac12
      - i \frac{\al k}{8\pi} \sigma(r)
\end{equation*}
acting on $\ltr = L^2(\R^+; dr)$, with  
\[
\sigma(r)=\frac{1-e^{-r^2/4}}{r^2/4},
\]
and domain
\begin{equation*}
    D(\mathscr{H}_{\al,k}) \eqdefa
\begin{cases}
\{ w\in\ltr :\ \pa_r^2 w,\ r^{-2}w,\ r^2 w \in\ltr \}, & |k|\ge 2, \\
\{ w\in\ltr :\ r^{1/2}\pa_r(r^{-3/2} w),\ r^{1/2}\pa_r(r^{-3/2} w),\ r^2 w \in\ltr \}, & |k|=1.
\end{cases}
\end{equation*}

As in Section~\ref{sect. 5}, we introduce the analytic family of type~(A)
\begin{equation}
    \mathscr{H}_{\al,k}^z \eqdefa
    z^{-2}\left(\pa_r^2 - (k^2 - \tfrac14) r^{-2}\right)
    - z^2 \frac{r^2}{16} + \frac12 - i \frac{\al k}{8\pi} \sigma(zr),
\end{equation}
with common domain $D(\mathscr{H}_{\al,k})$ for $z \in \cS = \{ re^{i\theta} : -\pi/8 < \theta < \pi/8,\ r>0 \}$.  
Then $\sigma_{\ltr}(\mathscr{H}_{\al,k})$ consists only of discrete eigenvalues and is independent of $z \in \cS$.

Choose  
\[
\Bar{z} = e^{i\cdot \sign(\al k)\cdot \frac{\pi}{16}}.
\]
A direct computation yields
\begin{equation}\label{eq A.15}
\begin{split}
-\Re \langle \mathscr{H}_{\al,k}^{\Bar{z}} w,\, w \rangle_{\ltr}
&= \cos\frac{\pi}{8}
  \Big( \|\pa_r w\|_{\ltr}^2 + (k^2-\tfrac14)\|r^{-1}w\|_{\ltr}^2 
        + \tfrac{1}{16}\|r w\|_{\ltr}^2 \Big) \\
&\quad - \frac12\|w\|_{\ltr}^2
      - \frac{|\al k|}{8\pi}
        \int_0^\infty \Im \sigma(r e^{i\frac{\pi}{16}})\,|w|^2\,dr \\
&\ge C_1 \int_0^\infty 
       \bigl(k^2 r^{-2} + r^2 + \al k\min\{r^2, r^{-2}\} - C_2\bigr)
       |w|^2\,dr \\
&\ge C_3\,(|\al k|^{1/2} - C_4)\,\|w\|_{\ltr}^2.
\end{split}
\end{equation}
Here in the first equality we used  
\[
\Im \sigma(re^{i\theta}) = \sign(\theta)\, \Im \sigma(re^{i|\theta|}),
\]
and the bound (Lemma~7.2 of \cite{LWZ})
\[
-\Im \sigma(r e^{i\theta}) \gtrsim \sin\theta \min\{ r^2, r^{-2} \},
\quad \text{\rm for}\  r>0,\ 0<\theta<\pi/4.
\]

Estimate \eqref{eq A.15}, together with $|\al|\gg 1$, implies  
\[
s_{\ltr}(\mathscr{H}_{\al,k}) \le -C |\al k|^{1/2},
\]
which completes the proof of Lemma~\ref{Lem A.4}.
\end{proof}

\renewcommand{\theequation}{\thesection.\arabic{equation}}
\setcounter{equation}{0}
\section{Estimates for the 3D linearized equations corresponding to \eqref{eq 1.14}}
This section is devoted to the decay and smoothing properties of the linearized equations associated with \eqref{eq 1.14}. Denote by $\oome(\tau)=\cS(\tau,s)\oome(s)$ the solution of 
\begin{subequations}\label{eq 3.1}
\begin{gather}
    \pa_\tau \begin{pmatrix}
        \Omega^r\\\Omega^\theta
    \end{pmatrix}
    -(e^\tau\pa_z^2+ \sL - \al \Pi)\begin{pmatrix}
        \Omega^r\\\Omega^\theta
    \end{pmatrix} -\al \begin{pmatrix}
        R^r(\oome)\\R^\theta(\oome)
    \end{pmatrix} =0,\\
    \pa_\tau \Omega^z -(e^\tau\pa_z^2+L-\al \Lambda)\Omega^z- \al R^z(\oome) = 0,\\
    \with \oome(\tau=s)=\oome(s).
\end{gather}
\end{subequations}
To investigate the dissipation in the $z$-variable, introduce the Fourier multipliers
\[
\cM_{\tau,s} \eqdefa \exp\big((e^{\frac{\tau}{2}}-e^{\frac{s}{2}})|D_z|\big),
\]
which commute with $\cS(\tau,s)$.
Recalling that $\UU = \bgrad \times (-\blap)^{-1} \oome$, one obtains the elliptic system
\begin{equation*}
\left\{
    \begin{aligned}
        &-(\blap-r^{-2})U^r = r^{-1}\pa_\theta\Omega^z - e^{\frac{\tau}{2}}\pa_z\Omega^\theta -2r^{-2}\pa_\theta U^\theta, \\
        &-(\blap-r^{-2})U^\theta = e^{\frac{\tau}{2}}\pa_z\Omega^r-\pa_r\Omega^z + 2r^{-2}\pa_\theta U^r,\\
        &-\blap U^z=r^{-1}\pa_r(r\Omega^\theta) - r^{-1}\pa_\theta \Omega^r,
    \end{aligned}
\right.
\end{equation*}
which, together with the definition of $\boldsymbol{R}$, implies
\begin{equation*}
\left\{
\begin{aligned}
&\text{if } \oome=\oome_{\neq}, \quad 
\text{then } \
   \UU=\UU_{\neq} \andf  \boldsymbol{R}(\oome)=\boldsymbol{R}_{\neq}(\oome_{\neq}), \\
&\text{if } \oome=\oome_{0}, \quad \
\text{then } \
   \UU=\UU_{0} \andf  \boldsymbol{R}(\oome)=\boldsymbol{R}_{0}(\oome_{0}),
\end{aligned}
\right.
\end{equation*}
and hence $B_z(X_0)$ and $B_z(X_{\neq})$ are invariant subspaces of $\cS(\tau,s)$. The following conservation property for $S(\tau,s)$ and for the nonlinear evolution \eqref{eq 1.14} will be needed. In particular, it will be shown that $B_z(X_0)^2\times Z_0$ and $B_z(X_{\neq})^2\times Z_{\neq}$ are invariant subspaces of $S(\tau,s)$ (and of the nonlinear evolution system \eqref{eq 1.14}).

\begin{Proposition}\label{Prop 3.1b}
 Let $\oome(\tau)\in C([0,T);B_z(Y))$ be a solution to either the linear equation \eqref{eq 3.1} or the nonlinear equation \eqref{eq 1.14}, with initial data $\oome(\tau=0)=\oome_0\in B_z(Y)$ satisfying $\Omega^z_0\in Z$, that is,
\begin{equation}\label{eq 3.3b}
\int_{\R^2}\int_{\T} \Omega_0^z \,d \xxi \,dz =0 \andf 
\int_{\R^2}\int_{\T} \xxi \Omega_0^z \,d \xxi \,dz =0 .
\end{equation}
Then $\Omega^z(\tau)\in Z$ for all $\tau\geq 0$. 
\end{Proposition}

\begin{proof}
Let $\Omega(\tau)$ be a sufficiently smooth solution of \eqref{eq 1.14}. In view of \eqref{eq 1.15} and \eqref{eq 1.16}, introduce
\begin{equation*}\begin{split}
    J_1 &\eqdefa \pa_\tau \Omega^z -(e^\tau\pa_z^2+L) \Omega^z,\\
    J_2 &\eqdefa \vg\cdot\gradxi \Omega^z + \UU^\sxxi \cdot \gradxi G- e^{\frac{\tau}{2}} G\pa_z U^z, \\
    J_3& \eqdefa \oome\cdot \bgrad U^z - \UU\cdot\bgrad \Omega^z.
\end{split}\end{equation*}
Using integration by parts and $\bdive\, \oome=\bdive\, \UU = \grad_\sxxi \cdot \vg =\pa_z G=0$, one obtains
\begin{gather*}
    \int_{\R^2}\int_{\T} J_1 \,d\xxi \,dz
    = \pa_\tau \int_{\R^2}\int_{\T} \Omega^z(\tau) \,d \xxi \,dz, \qquad  
    \int_{\R^2}\int_{\T} J_2 \,d\xxi \,dz
     =\int_{\R^2}\int_{\T} J_3 \,d\xxi \,dz =0,  \\
    \int_{\R^2}\int_{\T} \xxi J_1 \,d\xxi \,dz
    = \Big(\pa_\tau+\frac{1}{2}\Big) \int_{\R^2}\int_{\T} \xxi \Omega^z(\tau) \,d \xxi \,dz.
\end{gather*}
This yields the first identity in \eqref{eq 3.3b} for $\Omega^z(\tau)$.

For the second identity in \eqref{eq 3.3b}, it suffices to verify that
\begin{equation}\label{eq 3.4b}
    \int_{\R^2}\int_{\T} \xxi J_2 \,d\xxi \,dz
    =\int_{\R^2}\int_{\T} \xxi J_3 \,d\xxi \,dz=0.
\end{equation}
Indeed, integration by parts gives
\begin{equation*}
\begin{split}
    \int_{\R^2}\int_{\T} \xxi J_2 \,d\xxi \,dz
    &=  - \int_{\R^2}\int_{\T} \big(\vg\cdot\gradxi \xxi \big) \Omega^z \,d\xxi \,dz 
       -  \int_{\R^2}\int_{\T} \big(\UU^\sxxi\cdot\gradxi \xxi \big)G \,d\xxi \,dz \\
    &= - \int_{\R^2}\int_{\T} \vg \big(\gradxi^\perp \cdot \UU^\sxxi\big) \,d\xxi \,dz 
       - \int_{\R^2}\int_{\T} \UU^\sxxi  G \,d\xxi \,dz \\
    &= \int_{\R^2}\int_{\T} \UU^\sxxi \big(\gradxi^\perp \cdot\vg -G   \big) \,d\xxi \,dz =0,
\end{split}
\end{equation*}
and
\begin{equation*}
\begin{split}
   \int_{\R^2}\int_{\T} \xi_1 J_3 \,d\xxi \,dz
   &= \int_{\R^2}\int_{\T} U^1 \Omega^z - \Omega^1 U^z  \,d\xxi \,dz \\
   &= \int_{\R^2}\int_{\T} \left(U^1 (\pa_1 U^2-\pa_2 U^1) - (\pa_2 U^z - e^{\tau/2}\pa_z U^2) U^z \right)\,d\xxi \,dz\\
   &=  \int_{\R^2}\int_{\T} (U^1\pa_1 + e^{\tau/2}U^z \pa_z)U^2 \,d\xxi \,dz 
    = - \int_{\R^2}\int_{\T} U^2\pa_2 U^2 \,d\xxi \,dz =0.
\end{split}
\end{equation*}
As $\int_{\R^2}\int_{\T} \xi_2 J_3 \,d\xxi \,dz=0$ can be similarly deduced, we conclude \eqref{eq 3.4b} from the above identities. This completes the proof of Proposition~\ref{Prop 3.1b}.
\end{proof}

The main result of this section is stated below, and its proof will be given at the end of the section.

\begin{Proposition}\label{Prop 3.2}
 Let $1< p \leq 2\leq q <\oo$ and $0\leq s\leq\tau$. There exists a constant $C_\al>0$, depending only on $\al$, such that if $\oome(s)\in B_z(X_0)^2\times Z_0$, then
\begin{equation}\label{eq 3.3a}
\begin{split}
&a(\tau-s)^{\f1p-\f1q} \|\cM_{\tau,s}\cS(\tau,s)\oome(s)\|_{B_z(L^q(\oo))}
  + a(\tau-s)^{\f1p}\|\bgrad \cM_{\tau,s}\cS(\tau,s)\oome(s)\|_{B_z(L^2(\oo))}  \\
&\qquad \leq C_\al \|\oome(s)\|_{B_z(L^p(\oo))}(1+\tau-s)e^{-(\tau-s)/2}.
\end{split}
\end{equation}
If $\oome(s) \in B_z(X_{\neq})^2\times Z_{\neq}$, then
\begin{equation}\label{eq 3.4a}
\begin{split}
&a(\tau-s)^{\f1p-\f1q} \|\cM_{\tau,s}\cS(\tau,s)\oome(s)\|_{B_z(L^q(\oo))}
 + a(\tau-s)^{\f1p}\|\bgrad \cM_{\tau,s}\cS(\tau,s)\oome(s)\|_{B_z(L^2(\oo))}  \\
&\qquad \leq C_\al \|\oome(s)\|_{B_z(L^p(\oo))}e^{-\kappa |\al|^{1/2}(\tau-s)}.
\end{split}
\end{equation}
\end{Proposition}

\subsection{Estimates of $e^{\tau(\sL-\al \Pi)}$ and $e^{\tau(L-\al \La)}$}

In this subsection, the decay properties of $e^{\tau(\sL-\al \Pi)}$ and $e^{\tau(L-\al \La)}$ at different $\theta$-frequencies are investigated. Set
\begin{equation*}
\begin{split}
&\lap_k \eqdefa \lapxi\big|_{X_k}
= \pa_r^2+\frac{1}{r}\pa_r-\frac{k^2}{r^2},\qquad 
L_{k}\eqdefa L\big|_{X_k}
= \pa_r^2+\frac{1}{r}\pa_r-\frac{k^2}{r^2}+\frac{r}{2}\pa_r+1,\\
&\sL_k\eqdefa \sL\big|_{X_k^2}
= (L_k-r^{-2}){\rm Id}-\frac{2ik}{r^2} 
\begin{pmatrix}
0 & 1\\
-1 & 0 
\end{pmatrix},\qquad  
\Pi_k\eqdefa \Pi\big|_{X_k^2}
= \begin{pmatrix}
        i k S(r) & 0 \\
       -rS'(r) & i k S(r)
    \end{pmatrix},
\end{split}
\end{equation*}
with domains obtained by restricting
\begin{equation*}
\begin{split}
&D(L)\eqdefa \big\{f\in Y  : \  \lapxi f,\   \xxi\cdot\gradxi f \in Y  \big\},\\
&D(\sL)\eqdefa  \big\{(f_1,f_2)^\tr \in Y^2  : \  \lapxi f_i,\  \xxi\cdot\gradxi f_i,\  r^{-2}f_i,\  r^{-2}\pa_\theta f_i\in Y \quad \text{for } i=1,2  \big\}, \\
&D(\Pi)\eqdefa \big\{(f_1,f_2)^\tr \in Y^2  : \  S(r)\pa_\theta f_i,\ rS'(r)f_i \in Y \quad \text{for } i=1,2  \big\}.
\end{split}
\end{equation*}

\begin{Proposition}\label{Prop 3.3}
 There exist universal constants $\kappa,C>0$ and an $\al$-dependent constant $C_\al$ such that 
\begin{subequations}
\begin{gather}
    e^\tau \big\|e^{\tau(L-\al \La)}\big\|_{X_0\cap Y_0}
    + \big\|e^{\tau(L-\al \La)}\big\|_{X_0}  
    + C_\al^{-1} e^{\kappa |\al|^{1/2} \tau} 
      \big\|e^{\tau(L-\al \La)}\big\|_{X_{\neq} \cap Y_1}   \leq 1,  \label{eq 3.1a}\\
     \big\|e^{\tau(\sL-\al \Pi)}\big\|_{X_0^2} \leq C_\al (1+\tau)e^{-\tau/2}, \qquad 
      \big\|e^{\tau(\sL-\al \Pi)}\big\|_{X_{\neq}^2} \leq C_\al  e^{-\kappa |\al|^{1/2} \tau}. \label{eq 3.1b}
\end{gather}
\end{subequations}
\end{Proposition}

\begin{proof}
Since $\La_0=0$, Sturm–Liouville theory implies that $e^{-r^2/4}$ and $(r^2-4)e^{-r^2/4}$ are eigenfunctions of $L_0$, corresponding to the eigenvalues $0$ and $-1$, respectively. From this fact, the first two estimates in \eqref{eq 3.1a} follow. The last estimate in \eqref{eq 3.1a} and the bounds in \eqref{eq 3.1b} are direct consequences of Theorem~\ref{Thm 1.3}, Lemma~\ref{Lem 3.2}, and relatively compact perturbation theory. It remains to prove the first estimate in \eqref{eq 3.1b}. 

Observe that
\begin{equation}\label{eq 3.5}
    \sL_0-\al \Pi_0 = 
    \begin{pmatrix}
        L_1 & 0 \\ \al rS'(r) & L_1
    \end{pmatrix}.
\end{equation}
Let $f(\tau)=(f_1(\tau),f_2(\tau))^\tr
   =e^{\tau(\sL_0-\al\Pi_0)} f(0)$. By \eqref{eq 3.5} and Duhamel’s formula, one has
\begin{equation*}
\begin{split}
    f_1(\tau) &= e^{\tau L_1}f_1(0), \\
    f_2(\tau) 
    &= e^{\tau L_1}f_2(0)
       +\al \int_0^\tau e^{(\tau-s)L_1} rS'(r)e^{s L_1}f_1(0) \,ds .
\end{split}
\end{equation*}
By using $|rS'(r)|\lesssim 1$ and $\|e^{\tau L_1}\|_{X_0}\leq e^{-\tau/2}$, we get
\begin{equation*}
\begin{split}
     \|f_1(\tau)\|_{Y} &\leq e^{-\frac{\tau}{2}} \|f_1(0)\|_{Y}, \\
     \|f_2(\tau)\|_{Y}  
     &\leq e^{-\frac{\tau}{2}} \|f_2(0)\|_{Y}
        +C\al \int_0^\tau e^{-\frac{\tau-s}{2}}e^{-\frac{s}{2}} \,ds\, \|f_1(0)\|_{Y} \\
     &\leq e^{-\frac{\tau}{2}} \|f_2(0)\|_{Y} 
        +C\al\tau e^{-\frac{\tau}{2}} \|f_1(0)\|_{Y}.
\end{split}
\end{equation*}
Hence $\|e^{\tau(\sL_0-\al \Pi_0)}\|_{X_0^2} \lesssim_\al (1+\tau)e^{-\tau/2}$, which completes the proof of Proposition~\ref{Prop 3.3}.
\end{proof}

\subsection{Decay property of $S(\tau,s)$}

As a consequence of Proposition~\ref{Prop 3.3}, the following decay estimates hold.

\begin{Proposition}\label{Prop 3.4}
 There exist a universal constant $\kappa>0$ and an $\al$-dependent constant $C_\al>0$ such that for any $0\leq s\leq \tau$, 
\begin{subequations}
\begin{gather}
\|\cM_{\tau,s}\cS(\tau,s)\|_{B_z(X_0^2)\times Z_0} 
    \leq C_\al (1+\tau-s)e^{-(\tau-s)/2},\label{eq 3.8a}\\
\|\cM_{\tau,s}\cS(\tau,s)\|_{B_z(X_{\neq}^2)\times Z_{\neq}} 
    \leq C_\al e^{-\kappa |\al|^{1/2}(\tau-s)}.\label{eq 3.8b}
\end{gather}
\end{subequations}
\end{Proposition}

\begin{proof}
Let $\oome(\tau)=\cS(\tau,s)\oome(s)$, and define
\[
\boldsymbol{\omega}(\tau,\xxi,\zeta)
    =e^{e^\tau\zeta^2} \widehat{\oome}(\tau,\xxi,\zeta),
    \qquad
\boldsymbol{F}(\tau,\xxi,\zeta)
    =e^{e^\tau\zeta^2} \widehat{\boldsymbol{R}}(\tau,\xxi,\zeta).
\]
In view of \eqref{eq 3.1}, for each $\zeta\in\N$,
\begin{equation}\label{eq 3.12b}
\begin{split}
    \begin{pmatrix}
        \om^r \\ \om^\theta
    \end{pmatrix}(\tau) 
    &= e^{(\tau-s)(\sL-\al\Pi)}
       \begin{pmatrix}
         \om^r \\ \om^\theta
       \end{pmatrix}(s)
       + \int_s^\tau e^{(\tau-\rho)(\sL-\al\Pi)} 
       \begin{pmatrix}
        F^r \\ F^\theta
       \end{pmatrix} (\rho)\,d\rho, \\
   \om^z(\tau) 
   &= e^{(\tau-s)(L-\al\La)}\om^z(s) 
      + \int_s^\tau e^{(\tau-\rho)(L-\al\La)} F^z(\rho) \,d\rho.
\end{split}
\end{equation}
From Lemma~\ref{lem B.2} (with $\delta=\f12$), we deduce that 
\begin{equation}\label{eq 3.13b}
     \|\boldsymbol{F}(\cdot,\zeta)\|_Y
     \leq C_\al e^{\frac{\tau}{4}}|\zeta|^\f12 
     \|\boldsymbol{\om}(\cdot,\zeta)\|_Y.
\end{equation}

\noindent {\bf Step 1. Proof of \eqref{eq 3.8a}.}
Assume $\oome(s) \in B_z(X_0)^2\times Z_0$. Proposition~\ref{Prop 3.1b} then implies
\begin{equation*}
    \boldsymbol{\om}(\tau,\cdot,\zeta)\in 
    \begin{cases}
        X_0^2\times (X_0\cap Y_1), & \text{if } \zeta=0,\\
        X_0^2\times X_0, & \text{if } \zeta\neq0.
    \end{cases}
\end{equation*}
The proof is divided into the cases $\zeta=0$ and $\zeta\neq0$.

\noindent{\bf (1) Case $\zeta=0$.}
Proposition~\ref{Prop 3.3}, together with \eqref{eq 3.12b} and \eqref{eq 3.13b}, yields
\begin{equation*}
    \|\boldsymbol{\om}(\tau)\|_Y
    \lesssim_\al  (1+\tau-s)e^{-(\tau-s)/2}\|\boldsymbol{\om}(s)\|_Y,
\end{equation*}
and hence, by the definition of $\boldsymbol{\om}$,
\begin{equation}\label{eq 3.11}
    \|\widehat{\oome}(\tau,\cdot,\zeta)\|_Y 
    \lesssim_\al   (1+\tau-s)e^{-(\tau-s)/2} 
    \|\widehat{\oome}(s,\cdot,\zeta)\|_Y.
\end{equation}

\noindent{\bf (2) Case $\zeta\neq 0$.}
In this case, Proposition~\ref{Prop 3.3} together with \eqref{eq 3.12b} and \eqref{eq 3.13b} gives
\begin{equation*}
    \|\boldsymbol{\om}(\tau)\|_Y
    \lesssim_\al \|\boldsymbol{\om}(s)\|_Y 
    + \int_s^\tau e^{\frac{\rho}{4}}|\zeta|^\f12 
      \|\boldsymbol{\om}(\rho)\|_Y \,d\rho.
\end{equation*}
Applying Gronwall’s inequality yields
\begin{equation*}
\|\boldsymbol{\om}(\tau)\|_Y 
   \lesssim_\al \exp\big((e^{\tau/4}-e^{s/4})|\zeta|^{1/2}\big) 
   \|\boldsymbol{\om}(s)\|_Y.
\end{equation*}
Consequently, it comes out
\begin{equation}\label{eq 3.12}
\begin{split}
\|\widehat{\oome}(\tau,\cdot,\zeta)\|_Y 
&\lesssim_\al  
   e^{- (e^\tau-e^s)\zeta^2 + (e^{\tau/4}-e^{s/4})|\zeta|^{1/2}} 
   \|\widehat{\oome}(s,\cdot,\zeta)\|_Y\\
&\lesssim_\al e^{-(\tau-s)/2} 
   e^{-(e^{\tau/2}-e^{s/2})|\zeta|} 
   \|\widehat{\oome}(s,\cdot,\zeta)\|_Y,
\end{split}
\end{equation}
where we used, for $0\leq s\leq \tau$ and $|\zeta|\geq 1$,
\begin{equation}\label{eq 3.13}
        - (e^\tau-e^s)\zeta^2 
        + (e^{\tau/4}-e^{s/4})|\zeta|^{1/2}
        \leq C-(\tau-s)/2 -(e^{\tau/2}-e^{s/2})|\zeta|.
\end{equation}
Combining \eqref{eq 3.11} and \eqref{eq 3.12} leads to \eqref{eq 3.8a}.

\noindent {\bf Step 2. Proof of \eqref{eq 3.8b}.}   
Assume $\oome(s) \in B_z(X_{\neq})^2\times Z_{\neq}$. Then Proposition~\ref{Prop 3.1b} yields
\begin{equation*}
    \boldsymbol{\om}(\tau,\cdot,\zeta)\in 
    \begin{cases}
        X_{\neq}^2\times (X_{\neq}\cap Y_1), & \text{if } \zeta=0,\\
        X_{\neq}^2\times X_{\neq}, & \text{if } \zeta\neq0.
    \end{cases}
\end{equation*}
The argument proceeds as in {\bf Step~1}:

\noindent{\bf (1) Case $\zeta=0$.}
Using Proposition~\ref{Prop 3.3} to replace $(1+\tau-s)e^{-(\tau-s)/2}$ in \eqref{eq 3.11} by $e^{-\kappa |\al|^{1/2} (\tau-s)}$ gives the desired bound.

\noindent{\bf (2) Case $\zeta\neq 0$.}  
It suffices to replace \eqref{eq 3.13} by 
\begin{equation*}
        - (e^\tau-e^s)\zeta^2 + (e^{\tau/4}-e^{s/4})|\zeta|^{1/2}
        \leq C_\al-(e^{\tau/2}-e^{s/2})|\zeta| - \kappa |\al|^{1/2} (\tau-s).
\end{equation*}
The resulting estimates in the two cases yield \eqref{eq 3.8b} and complete the proof of Proposition~\ref{Prop 3.4}. 
\end{proof} 

\subsection{Smoothing property of $S(\tau,s)$}

Although the following lemma is standard, a brief proof is included for completeness.

\begin{Lemma}\label{Lem 3.3}
For $1\leq p \leq q \leq \oo$ and $\beta\in \N^2$, the estimate
\begin{equation}\label{eq 3.6}
    \|\gradxi^\beta e^{\tau L} f\|_{L^q(\oo)} 
    \lesssim a(\tau)^{-\left(\f1p-\f1q+\frac{|\beta|}{2}\right)} 
    \| f\|_{L^p(\oo)}
\end{equation}
holds for all $\tau\geq 0$.
\end{Lemma}

\begin{proof}
By virtue of Appendix~A in \cite{gallay:2001a}, we have
\begin{equation*}
    e^{\tau L}f(\xxi)
    =\frac{1}{4\pi a(\tau)} 
     \int_{\R^2} e^{-\frac{|\sxxi-\scalebox{0.6}{\eeta} e^{-\tau/2}|^2}{4a(\tau)}} 
     f(\eeta) \,d\eeta.
\end{equation*}
For $\beta=0$, set $\phi=G^{-\frac{1}{2}}e^{\tau L} f$ and $\psi=G^{-\frac{1}{2}}f$. Then
\begin{equation*}
\phi(\xxi) 
= \frac{1}{4\pi a(\tau)} 
  \int_{\R^2} K_\tau(\xxi,\eeta)\psi(\eeta)\,d\eeta,
\end{equation*}
with 
\begin{equation*}
\begin{split}
      K_\tau(\xxi,\eeta)
      &\eqdefa \exp\big(|\xxi|^2/8 -|\xxi-e^{-\tau/2}\eeta|^2/(4a(\tau)) -|\eeta|^2/8\big)\\
      &= \exp\Big(- \frac{(1+e^{-\tau})(|\xxi|^2+|\eeta|^2) -4e^{-\tau/2}\xxi\cdot\eeta}{8a(\tau)}  \Big) \\
      &= \exp\Big(-\frac{a^2(\tau/2)(|\xxi|^2 +|\eeta|^2) + 2e^{-\tau/2}|\xxi-\eeta|^2 }{8a(\tau)}\Big).  
\end{split}
\end{equation*}
Since $a(\tau/2) \sim a(\tau)$,
\begin{equation*}
\|\phi\|_{L^q_\sxxi} 
\lesssim a(\tau)^{-1} \|K_\tau\|_{ L^q_\sxxi (L^{p'}_{\scalebox{0.6}{\eeta}})} 
\|\psi\|_{L^p_{\scalebox{0.6}{\eeta}}} 
\lesssim a(\tau)^{-1} 
\big(a(\tau)/a^2(\tau/2) \big)^{-\frac{1}{p'}-\frac{1}{q}}
\|\psi\|_{L^p_{\scalebox{0.6}{\eeta}}} 
\lesssim a(\tau)^{\f1q-\f1p}\|\psi\|_{L^p_{\scalebox{0.6}{\eeta}}},
\end{equation*}
which is equivalent to \eqref{eq 3.6} when $\beta=0$. For general $\beta$, it suffices to take $\gradxi^\beta$ on the Gaussian kernel $e^{-\frac{|\sxxi-\scalebox{0.6}{\eeta} e^{-\tau/2}|^2}{4a(\tau)}}$, and then estimate the resulting kernel $K_{\beta,\tau}(\xxi,\eeta)$ in the same way. This completes the proof of Lemma~\ref{Lem 3.3}.
\end{proof}

\begin{Proposition}\label{Prop 3.5}
 Let $1< p\leq 2 \leq q <\oo$. There exist a constant $C_\al>0$ and a small parameter $\delta(\al)\ll 1$, depending only on $\al$, such that for any $s\leq \tau \leq s+ \delta(\al)$,
\begin{equation}\label{eq 3.15d}
\begin{split}
&(\tau-s)^{\f1p-\f1q} \|\cM_{\tau,s}\cS(\tau,s)\oome(s)\|_{B_z(L^q(\oo))}  \\
&\qquad + (\tau-s)^{\f1p}\|\bgrad \cM_{\tau,s}\cS(\tau,s)\oome(s)\|_{B_z(L^2(\oo))} 
\leq C_\al \|\oome(s)\|_{B_z(L^p(\oo))} .
\end{split}
\end{equation}
\end{Proposition}

\begin{proof}
Let $\oome(\tau)=\cS(\tau,s)\oome(s)$, and define
\[
\boldsymbol{\omega}(\tau,\xxi,\zeta)
   =e^{e^\tau\zeta^2} \widehat{\oome}(\tau,\xxi,\zeta),\qquad
\uu(\tau,\xxi,\zeta)
   =e^{e^\tau\zeta^2} \widehat{\UU}(\tau,\xxi,\zeta).
\]
For each $\zeta\in\N$, we find
\begin{equation}
    (\pa_\tau-L) \begin{pmatrix}
        \boldsymbol{\om}^\sxxi \\\om^z
    \end{pmatrix} 
    = - \al \vg\cdot \grad_\sxxi \begin{pmatrix}
         \boldsymbol{\om}^\sxxi \\ \om^z
    \end{pmatrix} 
      -i\al e^{\frac{\tau}{2}} \zeta G\begin{pmatrix}
        \uu^\sxxi \\ u^z
    \end{pmatrix} 
      + \al \begin{pmatrix}
         \boldsymbol{\om}^\sxxi  \cdot \grad_\sxxi \vg\\
        -\uu^\sxxi \cdot \grad_\sxxi G
    \end{pmatrix} 
    \eqdefa \boldsymbol{F}(\boldsymbol{\om}).
\end{equation}
Introduce the norm
\begin{equation}
\|\boldsymbol{\om}\|_{\mathcal{N}} \eqdefa \sup\limits_{s\leq\tau\leq s+\delta} 
\big( \| 
\boldsymbol{\om}\|_{L^p(\oo)} 
+ (\tau-s)^{\f1p-\f1q}\|\boldsymbol{\om}\|_{L^q(\oo)} 
+(\tau-s)^{\f1p}\|\grad_\sxxi \boldsymbol{\om}\|_{L^2(\oo)}\big).
\end{equation}
Lemma~\ref{lem B.1} implies, for $0<\tau-s\leq \delta\ll1$,
\begin{equation}\label{eq 3.15}
\begin{split}
    \|\boldsymbol{F}\|_{L^p(\oo)}
    &\lesssim_\al \|\vv^{\G}\|_{L^{2p/(2-p)}_\sxxi} 
                  \|\grad_\sxxi \boldsymbol{\om}\|_{L^2(\oo)} 
      + e^{\frac{\tau}{2}}|\zeta| \|\uu\|_{L^p_\sxxi} \|G\|_{L^\oo(\oo)}\\
    &\quad +\|\boldsymbol{\om}\|_{L^p(\oo)} \|\gradxi \vg\|_{L^\oo_{\sxxi}} 
      + \|\gradxi G\|_{L^{4p/(4-p)}(\oo)} \|\uu\|_{L^4_\sxxi}  \\
    &\lesssim_\al (\tau-s)^{-1/p} \|\boldsymbol{\om}\|_{\cN}.
\end{split}
\end{equation}
Using Lemma~\ref{Lem 3.3} and the fact that $a(\tau)\lesssim \tau$ for $\tau\leq 1$, we get
\begin{equation*}
\begin{split}
\|\boldsymbol{\om}\|_{\cN} 
&\leq \|e^{(\tau-s) L }\boldsymbol{\om}(s) \|_{\cN} 
     + \Big\|\int_s^\tau \boldsymbol{F}(\rho) \,d\rho  \Big\|_{\cN}\\
&\leq C \|\boldsymbol{\om}(s)\|_{L^p(\oo)} 
   + C_\al \|\boldsymbol{\om}\|_{\cN} 
     \sup_{0\leq \tau \leq \delta} 
     \Big(  \int_s^\tau (\rho-s)^{-\f1p} d\rho \\
&\quad + (\tau-s)^{\frac{1}{p}-\frac{1}{q}} \int_s^\tau (\tau-\rho)^{\frac{1}{q}-\frac{1}{p}} (\rho-s)^{-\f1p} d\rho 
      + (\tau-s)^{\f1p} \int_s^\tau (\tau-\rho)^{-\f1p}(\rho-s)^{-\f1p}d\rho \Big)\\
&\leq C \|\boldsymbol{\om}(s)\|_{L^p(\oo)} 
    + C_\al\delta^{1-\f1p}\|\boldsymbol{\om}\|_{\cN}.
\end{split}
\end{equation*}
Thus by choosing $\delta>0$ sufficiently small so that $C_\al\delta^{1-\f1p} \leq \f12$, it follows that
\[
\|\boldsymbol{\om}\|_{\cN} \leq  2C \|\boldsymbol{\om}(s)\|_{L^p(\oo)}.
\]

Finally, in view of the definition of $\boldsymbol{\om}$ and the bound, valid for $0\leq \tau-s \leq 1$ and all $k,\zeta\in\N$,
\begin{equation*}
\big(e^{\frac{\tau}{2}}|\zeta|\big)^k 
\exp\big(-(e^\tau-e^s)\zeta^2+ (e^{\frac{\tau}{2}}-e^{\frac{s}{2}})|\zeta|  \big) 
\leq C_k (\tau-s)^{-k/2},
\end{equation*}
we achieve \eqref{eq 3.15d} and thus complete the proof of Proposition~\ref{Prop 3.5}.
\end{proof}

The proof of Proposition~\ref{Prop 3.2} can now be given.

\begin{proof}[\bf Proof of Proposition \ref{Prop 3.2}]
Let $\delta=\delta(\al)\ll 1$ be as in Proposition~\ref{Prop 3.5}. If $\tau-s \leq  2\delta$, then Proposition~\ref{Prop 3.5} gives
\begin{equation}\label{eq 3.5a}
\begin{split}
&a(\tau-s)^{\f1p-\f1q} \|\cM_{\tau,s}\cS(\tau,s)\oome(s)\|_{B_z(L^q(\oo))}
 + a(\tau-s)^{\f1p}\|\bgrad \cM_{\tau,s}\cS(\tau,s)\oome(s)\|_{B_z(L^2(\oo))}  \\
&\leq C_\al \|\oome(s)\|_{B_z(L^p(\oo))} 
 \leq C_{\al}'M(\tau-s)\|\oome(s)\|_{B_z(L^p(\oo))},
\end{split}
\end{equation}
where $M(\tau)= (1+\tau)e^{-\tau/2}$ if $\oome\in B_z(X_0)^2\times Z_0$, and $M(\tau)= e^{-\kappa|\al|^{1/2}\tau}$ if $\oome \in B_z(X_{\neq})^2\times Z_{\neq}$.

If $\tau-s \geq 2\delta$, then Propositions~\ref{Prop 3.4} and \ref{Prop 3.5}, together with $1\geq a(\tau-s)\geq a(2\delta)\gtrsim_\al 1$, imply
\begin{equation}\label{eq 3.6a}
\begin{split}
&a(\tau-s)^{\f1p-\f1q} \|\cM_{\tau,s}\cS(\tau,s)\oome(s)\|_{B_z(L^q(\oo))}
 + a(\tau-s)^{\f1p}\|\bgrad \cM_{\tau,s}\cS(\tau,s)\oome(s)\|_{B_z(L^2(\oo))}  \\
&\leq C_\al \|\cM_{\tau-\delta/2,s}\cS(\tau-\delta/2,s)\oome(s)\|_{B_z(L^2(\oo))}    \\
&\leq C_\al' M(\tau-s-\delta) \|\cM_{s+\delta/2,s}\cS(s+\delta/2,s)\oome(s)\|_{B_z(L^2(\oo))}\\
&\leq C_{\al}''M(\tau-s)\|\oome(s)\|_{B_z(L^p(\oo))}.
\end{split}
\end{equation}

Combining \eqref{eq 3.5a} and \eqref{eq 3.6a}, we infer \eqref{eq 3.3a} and \eqref{eq 3.4a}, and thus complete the proof of Proposition~\ref{Prop 3.2}.
\end{proof}

\section{Proof of Theorem \ref{Thm 1}.}
In this section, the proof of Theorem \ref{Thm 1} is presented.  
Decompose $\oome=\oome_0+\oome_{\neq}$ and rewrite \eqref{eq 1.14} via the Duhamel formula as
\begin{subequations}
\begin{gather}
\oome(\tau)=\cS(\tau,0)\oome_0 - \int_0^\tau \cS(\tau,s)\boldsymbol{B}[\UU,\oome](s)\d s, \\
\oome_0(\tau)=\cS(\tau,0)\oome_0(0) - \int_0^\tau \cS(\tau,s)\boldsymbol{B}_0[\UU,\oome](s)\d s,\\
\oome_{\neq}(\tau)=\cS(\tau,0)\oome_{\neq}(0) - \int_0^\tau \cS(\tau,s)\boldsymbol{B}_{\neq}[\UU,\oome](s)\d s.
\end{gather}
\end{subequations}

Recall that $a(\tau)=1-e^{-\tau}$. Define
\begin{equation}\label{S4eq2}
\begin{split}
\|\oome\|_{N(T)}
&\eqdefa \sup_{0\le \tau\le T} (1+\tau)^{-1} e^{\tau/2}
\Big(
\|\cM_\tau \oome_0\|_{B_z(L^2(\oo))}
+ a^{1/4}(\tau)\|\cM_\tau \oome_0\|_{B_z(L^4(\oo))} \\
&\qquad\qquad\qquad
+ a^{1/2}(\tau)\|\cM_\tau \bgrad \oome_0\|_{B_z(L^2(\oo))}
\Big) \\
&\quad + \sup_{0\le \tau\le T} e^{\kappa|\al|^{1/2}\tau}
\Big(
\|\cM_\tau \oome_{\neq}\|_{B_z(L^2(\oo))}
+ a^{1/4}(\tau)\|\cM_\tau \oome_{\neq}\|_{B_z(L^4(\oo))} \\
&\qquad\qquad\qquad
+ a^{1/2}(\tau)\|\cM_\tau \bgrad \oome_{\neq}\|_{B_z(L^2(\oo))}
\Big).
\end{split}
\end{equation}

\begin{Lemma}\label{Lem 4.1}
There exists $C>0$ such that for all $0\le \tau\le T$,
\begin{subequations}
\begin{gather}
\|\cM_\tau \boldsymbol{B}_0[\UU,\oome]\|_{B_z(L^{4/3}(\oo))}
\le C\, a(\tau)^{-1/2} (1+\tau)^2 e^{-\tau} \|\oome\|_{N(T)}^2, \label{eq 4.3a} \\
\|\cM_\tau \boldsymbol{B}_{\neq}[\UU,\oome]\|_{B_z(L^{4/3}(\oo))}
\le C\, a(\tau)^{-1/2} (1+\tau) e^{-(1/2+\kappa|\al|^{1/2})\tau} \|\oome\|_{N(T)}^2. \label{eq 4.3b}
\end{gather}
\end{subequations}
\end{Lemma}

\begin{proof}
For $\frac1p+\frac1q=\frac1r$,
\begin{equation*}
\begin{split}
&\|\cM_\tau fg\|_{B_z(L^r(\oo))}
= \sum_{\zeta\in\Z} \exp\big((e^{\tau/2}-1)|\zeta|\big)
\|\widehat{fg}(\zeta,\cdot)\|_{L^r(\oo)} \\
&\leq \sum_{\zeta,\eta\in\Z} 
\exp\big((e^{\tau/2}-1)|\zeta-\eta|\big)
\exp\big((e^{\tau/2}-1)|\eta|\big)
\|\widehat{f}(\zeta-\eta,\cdot)\|_{L^p_\sxxi}
\|\widehat{g}(\eta,\cdot)\|_{L^q(\oo)} \\
&= \|\cM_\tau f\|_{B_z(L^p_\sxxi)}\,
   \|\cM_\tau g\|_{B_z(L^q(\oo))},
\end{split}
\end{equation*}
from which and  Lemma \ref{lem B.2}, we infer
\begin{equation}\label{eq 4.5}
\begin{split}
&\|\cM_\tau \boldsymbol{B}_0[\UU,\oome]\|_{B_z(L^{4/3}(\oo))}\\
&\le \|\cM_\tau (\UU\cdot \bgrad\oome)\|_{B_z(L^{4/3}(\oo))}
+ \|\cM_\tau (\oome\cdot \bgrad\UU)\|_{B_z(L^{4/3}(\oo))} \\
&\le \|\cM_\tau \UU\|_{B_z(L^4_\sxxi)}
\|\cM_\tau \bgrad\oome\|_{B_z(L^2(\oo))}
+ \|\cM_\tau \oome\|_{B_z(L^4(\oo))}
\|\cM_\tau \bgrad\UU\|_{B_z(L^2_\sxxi)} \\
&\le \|\cM_\tau \oome\|_{B_z(L^2(\oo))}
\|\cM_\tau \bgrad\oome\|_{B_z(L^2(\oo))}
+ \|\cM_\tau \oome\|_{B_z(L^4(\oo))}
\|\cM_\tau \oome\|_{B_z(L^2(\oo))}.
\end{split}
\end{equation}
Since $a(\tau)\le 1$, estimate \eqref{eq 4.3a} follows directly from \eqref{eq 4.5} and the definition of $N(T)$.

For \eqref{eq 4.3b}, note that
\[
\boldsymbol{B}_{\neq}[\UU,\oome]
= \boldsymbol{B}[\UU_0,\oome_{\neq}]
+ \boldsymbol{B}[\UU_{\neq},\oome_0]
+ \cP_{\neq}\boldsymbol{B}[\UU_{\neq},\oome_{\neq}],
\]
and the estimate follows from the same argument as above.  
This completes the proof of Lemma \ref{Lem 4.1}.
\end{proof}

Now the proof of Theorem \ref{Thm 1} can be given.

\begin{proof}[{\bf Proof of Theorem \ref{Thm 1}}]
Only the {\it a priori} estimates for sufficiently smooth solutions of \eqref{eq 1.14} are given.
Proposition \ref{Prop 3.2} yields
\[
\|\cS(\tau,0)\oome(0)\|_{N(T)}
\le C_\al \|\oome(0)\|_{B_z(Y)}.
\]
Using Proposition \ref{Prop 3.2} together with Lemma \ref{Lem 4.1}, one has
\begin{equation}\label{eq 4.6}
    \begin{split}
  &\Big\|\int_0^\tau \cS(\tau,s)\boldsymbol{B}_0[\UU,\oome](s) \d s\Big\|_{N(T)} \\
  &\lesssim_\al \|\oome\|_{N(T)}^2\sup_{0\leq \tau\leq T} \bigg(\frac{e^{\frac{\tau}{2}}}{1+\tau} \int_0^\tau (1+\tau-s)e^{-\frac{\tau-s}{2}} \times  a^{-\f12}(s)(1+s)^2e^{-s}   \\
  &\qquad \qquad \qquad \qquad  \times \Big( a^{-\f14} (\tau-s)  + a^\f14(\tau) a^{-\f12}(\tau-s)+ a^{\f12}(\tau)a^{-\f34}(\tau-s) \Big) \d s \bigg) \\
  &\lesssim_\al \|\oome\|_{N(T)}^2\sup_{0\leq \tau\leq T} \int_0^\tau e^{-\frac{s}{4}}   a^{-\f12} (s)\\
  &\qquad \qquad \qquad  \qquad\qquad  \times \Big( a^{-\f14} (\tau-s) + a^\f14(\tau) a^{-\f12}(\tau-s)+ a^{\f12}(\tau)a^{-\f34}(\tau-s) \Big) \d s \\
  &\lesssim_\al  a^{1/4}(T)\|\oome\|_{N(T)}^2,
    \end{split}
\end{equation}
and 
\begin{equation}\label{eq 4.7a}
    \begin{split}
  &\Big\|\int_0^\tau \cS(\tau,s)\boldsymbol{B}_{\neq}[\UU,\oome](s) \d s\Big\|_{N(T)} \\
  &\lesssim_\al \|\oome\|_{N(T)}^2\sup_{0\leq \tau\leq T} \bigg(e^{\kappa|\al|^{1/2}\tau} \int_0^\tau e^{-\kappa|\al|^{1/2}(\tau-s)}\times a^{-\f12}(s) (1+s)e^{-(\f12+\kappa|\al|^{1/2})s}   \\
  &\qquad \qquad \qquad \qquad \qquad \times \Big( a^{-\f14}  (\tau-s) + a^\f14 (\tau)a^{-\f12}(\tau-s)
  + a^{\f12}(\tau)a^{-\f34} (\tau-s)\Big) \d s \bigg) \\
    &\lesssim_\al \|\oome\|_{N(T)}^2\sup_{0\leq \tau\leq T} \int_0^\tau e^{-\frac{s}{4}}   a^{-\f12}(s) \\
  &\qquad \qquad \qquad  \qquad\qquad  \times \Big( a^{-\f14} (\tau-s) + a^\f14(\tau) a^{-\f12}(\tau-s)
  + a^{\f12}(\tau)a^{-\f34}(\tau-s) \Big) \d s \\
  &\lesssim_\al  a^{1/4}(T)\|\oome\|_{N(T)}^2.
    \end{split}
\end{equation}

Combining \eqref{eq 4.6}–\eqref{eq 4.7a} gives
\begin{equation}\label{eq 4.8}
\begin{split}
\|\oome\|_{N(T)}
\le C_\al \|\oome(0)\|_{B_z(Y)}
+ C_\al' a(T)^{1/4} \|\oome\|_{N(T)}^2
\le C_\al \|\oome(0)\|_{B_z(Y)}
+ C_\al' \|\oome\|_{N(T)}^2.
\end{split}
\end{equation}

Thus, if $\|\oome(0)\|_{B_z(Y)}$ is sufficiently small, then
\[
\|\oome\|_{N(T)} \lesssim_\al \|\oome(0)\|_{B_z(Y)}, \qquad \forall\, T<\infty,
\]
which completes the proof of Theorem \ref{Thm 1}.
\end{proof}

\appendix
\section{Basic calculus of complex deformations}\label{basic calculus of complex deformation}
Part of this Appendix is adapted from Section~7 of \cite{LWZ}. We begin by introducing the notations used in \cite{LWZ}:
\[
F_0(z) = e^z - z - 1, \quad  
F_1(z) = \frac{1-e^{-z}}{z}, \quad 
F_2(z) = e^{-z/2}, \quad
F_3(z) = \Bigl(\frac{2z^2}{F_0(z)} - 3 + 2z\Bigr)\frac{z}{F_0(z)}.
\]
It is straightforward to verify that $F_0, F_1,$ and $F_2$ are holomorphic in $\mathbb{C}$, whereas $F_3$ is holomorphic in
\[
\mathcal{S}_3 \eqdefa \big\{ r e^{i\theta} \,:\, -\pi/4 \leq \theta \leq \pi/4,\, r > 0 \big\},
\]
and meromorphic in $\mathbb{C}$. Moreover, in view of \eqref{S5eq1} and \eqref{S5eq2}, we have
\[
f(r) = F_3(r^2/4), \qquad \sigma(r) = F_1(r^2/4), \qquad g(r) = F_2(r^2/4), 
\quad \text{for } r \in \mathbb{R}^+,
\]
so that $F_3(z^2/4)$, $F_1(z^2/4)$, and $F_2(z^2/4)$ provide analytic extensions of $f(r)$, $\sigma(r)$, and $g(r)$ to the domain $\mathcal{S}_3$.

\begin{Lemma}\label{Lem A.1}
 For any $\zeta \in \mathcal{S}= \{ r e^{i\theta} : -\pi/8 < \theta < \pi/8, \ r > 0 \}$, there hold
\begin{subequations}
    \begin{gather}
        \sup_{r \in \mathbb{R}^+} \bigl|\zeta^2 f(\zeta r) - 8r^{-2}\bigr| \;\lesssim\; |\zeta|^2,  \label{estimate of f,ineq}\\
        \sup_{r \in \mathbb{R}^+} \bigl|\sigma(\zeta r) - 1 + \zeta^2 r^2/8\bigr|
        \;\lesssim\; \min\{|\zeta|^4 r^4,\; |\zeta|^2 r^2\}. \label{estimate of sig,ineq}
    \end{gather}
\end{subequations}
\end{Lemma}
\begin{proof}
We write $z = \zeta^2 r^2/4 \in \mathcal{S}_3$, and then observe that
\[
\zeta^2 f(\zeta r) - 8r^{-2} = \zeta^2\bigl(F_3(z) - 2/z\bigr),
\qquad 
\sigma(\zeta r) - 1 + \zeta^2 r^2/8 = F_1(z) - 1 + z/2.
\]
Thus it suffices to prove that
\[
 |F_3(z) - 2/z| \lesssim 1,
\qquad 
 |F_1(z) - 1 + z/2| \lesssim \min\{|z|^2,\, |z|\},
\quad \text{for } z \in \mathcal{S}_3.
\]
Since $z=0$ is a removable singularity of both $F_3(z) - 2/z$ and $F_1(z)$,  
we get by Taylor expansion that 
\begin{equation}\label{F_3 near 0}
|F_3(z) - 2/z| \lesssim 1,
\qquad 
|F_1(z) - 1 + z/2| \lesssim |z|^2,
\quad \text{for } |z| \leq 4.
\end{equation}
For $|z| \geq 4$, one has
\begin{equation*}
\begin{aligned}
|F_0(z)| 
&\geq |e^z| - |z| - 1 
 \;\geq\; e^{|z|/2} - |z| - 1 \andf
|e^{-z}|
\leq e^{-|z|/2},
\end{aligned}
\end{equation*}
which together with the definitions of the functions $F_i$, implies
\begin{equation}\label{F_3 near +oo}
 |F_3(z) - 2/z| \lesssim 1,
\qquad 
 |F_1(z) - 1 + z/2| \lesssim |z|,
\quad \text{for } |z| \geq 4.
\end{equation}
Combining \eqref{F_3 near 0} and \eqref{F_3 near +oo} gives the desired bounds
\eqref{estimate of f,ineq} and \eqref{estimate of sig,ineq}.
\end{proof}

\begin{Lemma}\label{Lem A.2}
 For any $\zeta \in \mathcal{S}_4 \eqdefa \{ r e^{i\theta} : \pi/16 \leq \theta \leq \pi/8,\, r > 0 \}$, 
there holds
\begin{equation}\label{estimate of coercivity of potential, ineq}
    \Im\big(\zeta^2 (1-\sigma(\zeta r))\big) 
    \gtrsim \min\big\{|\zeta|^4 r^2,\; |\zeta|^2\big\}.
\end{equation}
\end{Lemma}
\begin{proof}
We write $z=\zeta^2 r^2/4 = x + iy$, so that $y = cx$ with 
$c \in [\tan(\pi/8),\, 1]$, and observe that
\[
\zeta^2 (1-\sigma(\zeta r)) = 4r^{-2} z(1-F_1(z)),
\]
from which we obtain
\begin{equation}\label{eq A.5}
    \begin{split}
\Im\big(\zeta^2 (1-\sigma(\zeta r))\big)
&= 4r^{-2}\Im(e^{-z} - 1 + z)  \\
&= 4r^{-2}\bigl(cx - e^{-x}\sin(cx)\bigr)
     \eqdefa 4r^{-2} h(x).
    \end{split}
\end{equation}
Note that $h(x) \geq cx(1-e^{-x}) > 0$ for $x>0$.  
Using Taylor expansion, we infer
\begin{equation}\label{eq A.6}
    h(x) \sim \mathcal{O}(x^2) \quad \text{for } |x| \leq 1,
    \qquad
    h(x) \sim \mathcal{O}(x) \quad \text{for } |x| \geq 1.
\end{equation}
Combining \eqref{eq A.5}, \eqref{eq A.6}, and the relation 
$|x| \thicksim |\zeta|^2 r^2$, we deduce that
\[
\Im\big(\zeta^2 (1-\sigma(\zeta r))\big)
   \gtrsim \min\{\, |\zeta|^4 r^2,\; |\zeta|^2 \,\}.
\]
This completes the proof of Lemma~\ref{Lem A.2}.
\end{proof}

\begin{Lemma}\label{Lem A.3}
 For any $\zeta \in \mathcal{S}_4 = \{ r e^{i\theta} : \pi/16 \leq \theta \leq \pi/8,\, r > 0 \}$ and any $|k|\geq 2$, one has
\begin{equation}\label{ineq, coercivity with nonlocal term}
    \begin{split}
&\int_0^\infty \Im\!\big(\zeta^2 (1-\sigma(\zeta r))\big)\, |f(r)|^2\, dr
   + \int_0^\infty \int_0^\infty
     K_k(r,s)\, \Im\!\big(\zeta^4 g(\zeta r) g(\zeta s)\big)
     f(r)\,\overline{f(s)}\, dr ds \\
&\quad \gtrsim \int_0^\infty \min\{\, |\zeta|^4 r^2,\; |\zeta|^2 \,\} \, |f(r)|^2\, dr.
    \end{split}
\end{equation}
\end{Lemma}
\begin{proof}
We set $\zeta=|\zeta|e^{i\theta}$. Then the left-hand side of 
\eqref{ineq, coercivity with nonlocal term} can be written as
\begin{equation*}
    \begin{split}
LHS &= \int_0^\infty \frac{4}{r^2}
\Bigl(\frac{|\zeta|^2\sin(2\theta)}{4}r^2 
      - e^{-\frac{|\zeta|^2\cos(2\theta)}{4}r^2}
        \sin\!\Bigl(\frac{|\zeta|^2\sin(2\theta)}{4}r^2\Bigr)\Bigr)
|f(r)|^2\,dr \\
&\quad + |\zeta|^4 \int_0^\infty \int_0^\infty 
      K_k(r,s)\, e^{-\frac{|\zeta|^2\cos(2\theta)}{8}(r^2+s^2)}
      \sin\!\Bigl(4\theta - \frac{|\zeta|^2\sin(2\theta)}{8}(r^2+s^2)\Bigr)
      f(r)\overline{f(s)}\,drds \\
&\eqdefa J_1 + J_2 - J_3,
    \end{split}
\end{equation*}
where we use the trigonometric identity (cf.\ \cite{LWZ})
\[
\sin(a-b-c)\sin a
= \sin(a-b)\sin(a-c) - \sin b\, \sin c,
\]
and introduce
\begin{equation}
    \begin{split}
J_1 &\eqdefa \int_0^\infty \frac{4}{r^2}
\Bigl(\frac{|\zeta|^2\sin(2\theta)}{4}r^2 
      - e^{-\frac{|\zeta|^2\cos(2\theta)}{4}r^2}
        \sin\Bigl(\frac{|\zeta|^2\sin(2\theta)}{4}r^2\Bigr)\Bigr)
|f(r)|^2\,dr, \\
J_2 &\eqdefa \frac{|\zeta|^4}{\sin(4\theta)}
\int_0^\infty \int_0^\infty K_k(r,s)\,
e^{-\frac{|\zeta|^2\cos(2\theta)}{8}(r^2+s^2)}
\sin\Bigl(4\theta-\frac{|\zeta|^2\sin(2\theta)}{8}r^2\Bigr) \\
&\qquad\qquad\qquad \times
\sin\Bigl(4\theta-\frac{|\zeta|^2\sin(2\theta)}{8}s^2\Bigr)
f(r)\overline{f(s)}\,drds, \\
J_3 &\eqdefa \frac{|\zeta|^4}{\sin(4\theta)}
\int_0^\infty \int_0^\infty K_k(r,s)\,
e^{-\frac{|\zeta|^2\cos(2\theta)}{8}(r^2+s^2)}
\sin\Bigl(\frac{|\zeta|^2\sin(2\theta)}{8}r^2\Bigr)
\\
&\qquad\qquad\qquad \times\sin\Bigl(\frac{|\zeta|^2\sin(2\theta)}{8}s^2\Bigr)
f(r)\overline{f(s)}\,drds.
    \end{split}
\end{equation}

Since $K_k >0$ and $\sin(4\theta)>0$ for $\zeta\in\mathcal{S}_4$, we immediately have $J_2\ge 0$.  
Applying the Cauchy inequality and symmetry yields
\begin{equation}\label{J_3}
\begin{split}
J_3 &\le \frac{|\zeta|^4}{\sin(4\theta)}
\int_0^\infty \int_0^\infty K_k(r,s)\,
e^{-\frac{|\zeta|^2\cos(2\theta)}{8}(r^2+s^2)} \\
&\quad \qquad \qquad \times 
\Bigl|\sin\!\Bigl(\tfrac{|\zeta|^2\sin(2\theta)}{8}r^2\Bigr)\Bigr|
\Bigl|\sin\!\Bigl(\tfrac{|\zeta|^2\sin(2\theta)}{8}s^2\Bigr)\Bigr|
\frac{|f(r)|^2 + |f(s)|^2}{2}\,drds \\
&= \frac{|\zeta|^4}{\sin(4\theta)}
\int_0^\infty \varpi(r)\,
e^{-\frac{|\zeta|^2\cos(2\theta)}{8}r^2}
\Bigl|\sin\!\Bigl(\tfrac{|\zeta|^2\sin(2\theta)}{8}r^2\Bigr)\Bigr|
|f(r)|^2\,dr,
\end{split}
\end{equation}
where
\begin{equation}\label{Upsilon}
\begin{split}
\varpi(r)
&\eqdefa \int_0^\infty K_k(r,s)\,
e^{-\frac{|\zeta|^2\cos(2\theta)}{8}s^2}
\Bigl|\sin\!\Bigl(\tfrac{|\zeta|^2\sin(2\theta)}{8}s^2\Bigr)\Bigr|\,ds \\
&\le \frac{|\zeta|^2\sin(2\theta)}{32}
\int_0^\infty
\min\Bigl\{\frac{r}{s},\frac{s}{r}\Bigr\}^2
(rs)^{1/2}
e^{-\frac{|\zeta|^2\cos(2\theta)}{8}s^2}
s^2\,ds \\
&= \frac{|\zeta|^2\sin(2\theta)}{32}
\Bigl(
r^{-3/2} \int_0^r s^{9/2} e^{-\frac{|\zeta|^2\cos(2\theta)}{8}s^2}\,ds
+ r^{5/2} \int_r^\infty s^{1/2} e^{-\frac{|\zeta|^2\cos(2\theta)}{8}s^2}\,ds
\Bigr) \\
&\le \frac{|\zeta|^2\sin(2\theta)}{32}
\Bigl(
\int_0^r s^3 e^{-\frac{|\zeta|^2\cos(2\theta)}{8}s^2}\,ds
+ r^2 \int_r^\infty s e^{-\frac{|\zeta|^2\cos(2\theta)}{8}s^2}\,ds
\Bigr) \\
&= \frac{\tan(2\theta)}{|\zeta|^2\cos(2\theta)}
\Bigl(1 - e^{-\frac{|\zeta|^2\cos(2\theta)}{8}r^2}\Bigr).
\end{split}
\end{equation}
Here, in the second line, we have applied the estimates $K_k(r,s)\le K_2(r,s)$ and $|\sin x|\le |x|$.

Inserting \eqref{Upsilon} into \eqref{J_3} gives
\begin{equation}
    J_3 \le 
    \frac{|\zeta|^4 \tan^2(2\theta)}{8\sin(4\theta)}
    \int_0^\infty r^2
    \Bigl(1 - e^{-\frac{|\zeta|^2\cos(2\theta)}{8}r^2}\Bigr)
    e^{-\frac{|\zeta|^2\cos(2\theta)}{8}r^2}
    |f(r)|^2\,dr.
\end{equation}

Thus, to establish \eqref{ineq, coercivity with nonlocal term}, it remains to show that for $\theta\in[\pi/16,\pi/8]$,
\begin{equation}
\begin{split}
\varrho(r)
&\eqdefa 
\frac{4}{r^2}
\Bigl(\frac{|\zeta|^2\sin(2\theta)}{4}r^2
      - e^{-\frac{|\zeta|^2\cos(2\theta)}{4}r^2}
        \sin\!\Bigl(\frac{|\zeta|^2\sin(2\theta)}{4}r^2\Bigr)\Bigr) \\
&\quad
- \frac{|\zeta|^4\tan^2(2\theta)}{8\sin(4\theta)}
  r^2\Bigl(1-e^{-\frac{|\zeta|^2\cos(2\theta)}{8}r^2}\Bigr)
   e^{-\frac{|\zeta|^2\cos(2\theta)}{8}r^2} \\
&\gtrsim \min\{|\zeta|^4 r^2,\; |\zeta|^2\}.
\end{split}
\end{equation}
Set $\Lambda = \frac{|\zeta|^2\cos(2\theta)}{4} r^2$.  
Using $|\sin x|\le |x|$, we obtain
\begin{equation*}
\begin{split}
\varrho(\Lambda)
&= |\zeta|^2 \cos(2\theta)
\Bigl(
\tan(2\theta)
- e^{-\Lambda}\frac{|\sin(\tan(2\theta)\Lambda)|}{\Lambda}
- \frac{\tan^2(2\theta)}{2\sin(4\theta)\cos^2(2\theta)}\,
\Lambda\bigl(e^{-\Lambda/2}-e^{-\Lambda}\bigr)
\Bigr)
\\
&\ge |\zeta|^2 \sin(2\theta)
\Bigl(
1 - e^{-\Lambda}
- \frac{\tan(2\theta)}{2\sin(4\theta)\cos^2(2\theta)}\,
\Lambda\bigl(e^{-\Lambda/2}-e^{-\Lambda}\bigr)
\Bigr)
\\
&= |\zeta|^2 \sin(2\theta)
\Bigl(
1 - e^{-\Lambda}
- \frac{1}{4\cos^4(2\theta)}\,
\Lambda\bigl(e^{-\Lambda/2}-e^{-\Lambda}\bigr)
\Bigr).
\end{split}
\end{equation*}
Since $\theta\in[\pi/16,\pi/8]$, it follows that
\[
\varrho(\Lambda)
\ge |\zeta|^2 \sin(\pi/8)
\bigl(1 - e^{-\Lambda} - \Lambda(e^{-\Lambda/2}-e^{-\Lambda})\bigr)
\eqdefa |\zeta|^2\sin(\pi/8)\,f(\Lambda).
\]
Noting that
\[
(e^\Lambda f(\Lambda))'
= \bigl(e^{\Lambda/2}-1-\Lambda/2\bigr)e^{\Lambda/2}+1
\ge 1,
\]
we have $f(\Lambda)> e^{-\Lambda} f(0)=0$, which together with Taylor expansion yields
\[
\varrho(\Lambda)
\ge |\zeta|^2 \sin(\pi/8)\, f(\Lambda)
\gtrsim |\zeta|^2 \min\{\Lambda,\,1\}
\gtrsim \min\{|\zeta|^4 r^2,\; |\zeta|^2\}.
\]
This completes the proof of Lemma~\ref{Lem A.3}.
\end{proof}

\renewcommand{\theequation}{\thesection.\arabic{equation}}
\setcounter{equation}{0}
\section{Basic estimate in self-similar coordinates}

In this Appendix, we collect several basic estimates in self-similar coordinates.  
These bounds are essentially adapted from Lemma~3.3 and Corollary~3.4 in \cite{BGH23}.  
We begin by recalling that
\begin{equation*}
\begin{split}
&\oome = \oome^\sxxi + \Omega^z \mathbf{e}_z,\qquad   
\UU \eqdefa \bgrad \times (-\blap)^{-1}\oome = \UU^\sxxi + U^z \mathbf{e}_z,
\end{split}
\end{equation*}
and  
\[
\boldsymbol{R}(\oome) \eqdefa R^r \mathbf{e}_r + R^\theta \mathbf{e}_\theta + R^z \mathbf{e}_z,
\]
with
\begin{equation*}
\begin{split}
&R^r(\oome)\eqdefa G e^{\frac{\tau}{2}}\pa_z U^r,\qquad  
R^\theta(\oome)\eqdefa G e^{\frac{\tau}{2}}\pa_z U^\theta, \\
&R^z(\oome)\eqdefa G e^{\frac{\tau}{2}}\pa_z U^z
   - \gradxi G\cdot\big(\UU^\sxxi - \gradxi^\perp \lapxi^{-1}\Omega^z\big).
\end{split}
\end{equation*}

In self-similar Cartesian coordinates, the Biot–Savart law becomes
\begin{equation}
\UU^\sxxi
    = e^{\frac{\tau}{2}}\pa_z(-\blap)^{-1}\oome^{\sxxi,\perp}
      - \gradxi^\perp (-\blap)^{-1}\Omega^z,
\qquad 
U^z = \gradxi^\perp \cdot (-\blap)^{-1}\oome^\sxxi.
\end{equation}

Taking the Fourier transform in the $z$–variable yields
\begin{equation}\label{B.S.law in cartesian-Fourier}
\begin{cases}
\widehat{\UU}^\sxxi
    = i e^{\frac{\tau}{2}}\zeta (e^{\tau}\zeta^2 - \lapxi)^{-1}
      \widehat{\oome}^{\sxxi,\perp}
      - \gradxi^\perp (e^{\tau}\zeta^2 - \lapxi)^{-1}\widehat{\Omega}^z, \\
\widehat{U}^z
    = \gradxi^\perp \cdot 
      (e^{\tau}\zeta^2 - \lapxi)^{-1}\widehat{\oome}^\sxxi.
\end{cases}
\end{equation}

\begin{Lemma}\label{lem B.1}
 Let $m>1$ and $\la>0$.
\begin{itemize}
\item[(i)] For $1< r\le\infty$ and $0<\delta\le \min\{\tfrac12,\,1-\tfrac1r\}$,
\begin{equation}\label{C.2}
    \|(\la^2-\lapxi)^{-1}f\|_{L^r_\sxxi}
    \lesssim \la^{-(\frac{2}{r}+2\delta)}\|f\|_{L^2(m)}.
\end{equation}

\item[(ii)] For $1< r\le 2$ and $0<\delta\le1-\tfrac1r$,
\begin{equation}\label{C.3}
    \|\gradxi(\la^2-\lapxi)^{-1}f\|_{L^r_\sxxi}
    \lesssim \la^{1-(\frac{2}{r}+2\delta)}\|f\|_{L^2(m)}.
\end{equation}

\item[(iii)] For $2< r<\infty$,
\begin{equation}\label{C.4}
\|\gradxi(-\lapxi)^{-1}f\|_{L^r_\sxxi}
    + \|\gradxi(\la^2-\lapxi)^{-1}f\|_{L^r_\sxxi}
    \lesssim \|f\|_{L^2(m)}.
\end{equation}

\item[(iv)] For $1< r\le 2$,
\begin{equation}\label{C.5}
\|\gradxi^2(-\lapxi)^{-1}f\|_{L^r_\sxxi}
    + \|\gradxi^2(\la^2-\lapxi)^{-1}f\|_{L^r_\sxxi}
    \lesssim \|f\|_{L^2(m)}.
\end{equation}
\end{itemize}
\end{Lemma}

\begin{proof}
Estimates \eqref{C.2}, \eqref{C.3}, and the first bound in \eqref{C.4} were obtained in Lemma~3.3 of \cite{BGH23}.  
We therefore only verify the remaining estimates.

By H\"older’s inequality, we have for $1\le q\le 2$,
\[
\|f\|_{L^q_\sxxi} \lesssim \|f\|_{L^2(m)}.
\]
Since  
\[
(\la^2-\lapxi)^{-1}f(\xxi)
    = \int_{\R^2} K(\la(\xxi-\boldsymbol{\eta})) f(\boldsymbol{\eta})\, d\boldsymbol{\eta},
\]
where $K$ is the kernel of $(1-\lapxi)^{-1}$, we write
\[
\gradxi(\la^2-\lapxi)^{-1} f
    = \gradxi(-\lapxi)^{-1}f
      - \la^2 (\la^2-\lapxi)^{-1}\gradxi(-\lapxi)^{-1}f.
\]
Young’s inequality and the first estimate in \eqref{C.4} then imply
\[
\|\gradxi(\la^2-\lapxi)^{-1}f\|_{L^r_\sxxi}
    \lesssim (1+\|K\|_{L^1_\sxxi})
             \|\gradxi(-\lapxi)^{-1}f\|_{L^r_\sxxi}
    \lesssim \|f\|_{L^2(m)}.
\]

The bound in \eqref{C.5} follows from the $L^p$–boundedness of the Riesz transform and the same argument used in part (iii).  
This completes the proof of Lemma~\ref{lem B.1}.
\end{proof}

\begin{Lemma}\label{lem B.2}
 Let $m>1$. Then there hold
\begin{itemize}
\item[(i)] For $2<p<\infty$ and $1<q\le 2$,
\begin{equation}\label{estimate of U}
\|\widehat{\UU}\|_{L^p_\sxxi}
    + \|(\gradxi,\, i e^{\frac{\tau}{2}}\zeta)\widehat{\UU}\|_{L^q_\sxxi}
    \lesssim \|\widehat{\oome}\|_{L^2(m)}.
\end{equation}

\item[(ii)] For any $0\le \delta<1$, there exists $C=C(\al,\delta)>0$, independent of $\zeta$, such that
\begin{equation}\label{estimate of R}
\|\widehat{\boldsymbol{R}}(\cdot,\zeta)\|_Y
    \le C\,(e^{\frac{\tau}{2}}|\zeta|)^{\delta}
       \|\widehat{\oome}(\cdot,\zeta)\|_Y.
\end{equation}
\end{itemize}
\end{Lemma}

\begin{proof}
Setting $\la=e^{\frac{\tau}{2}}\zeta$ in \eqref{C.2}–\eqref{C.5}, estimate \eqref{estimate of U} follows directly from \eqref{B.S.law in cartesian-Fourier} and \eqref{definition Bz}.  
Estimate \eqref{estimate of R} follows from the rapid decay of $G$ and $\gradxi G$, together with the same argument used in the proof of Corollary~3.4 in \cite{BGH23}, with $Y$ playing the role of $L^2(m)$.
\end{proof}

\section*{Acknowledgement}

 T. Li is partially supported by  National Natural Science Foundation of China under Grant 12421001.  P. Zhang is partially  supported by National Key R$\&$D Program of China under grant 2021YFa1000800 and by National Natural Science Foundation of China under Grant 12421001, 12494542 and 12288201.

\section*{Declarations}

\subsection*{Conflict of interest} The authors declare that there are no conflicts of interest.

\subsection*{Data availability}
This article has no associated data.

\end{document}